\pdfoutput=1
\documentclass[11pt,onecolumn,journal]{article}
\usepackage[bottom=2.5cm,top=2.5cm,left=2.5cm,right=2.5cm]{geometry}

\usepackage{amsmath,amssymb,amsthm,graphicx,caption,cases}	
\usepackage[bottom]{footmisc}	
\usepackage{algorithmic}

\usepackage{abstract,amsfonts,amsbsy,amssymb,amstext,mathrsfs,latexsym,color,url,amsthm,xcolor,authblk,fontenc,bm}

\definecolor{NUSBlue}{RGB}{0,61,124} 
\definecolor{NUSOrange}{RGB}{239,124,0}
\usepackage[colorlinks=false,allbordercolors=NUSOrange]{hyperref}

\DeclareOldFontCommand{\bf}{\normalfont\bfseries}{\mathbf}

\newcommand{\bI}{\mathbf{I}}

\newcommand{\bp}{\mathbf{p}}
\newcommand{\bz}{\mathbf{z}}
\newcommand{\br}{\mathbf{r}}
\newcommand{\bv}{\mathbf{v}}

\newcommand{\bg}{\mathbf{g}}

\newcommand{\bu}{\mathbf{u}}

\newcommand{\bx}{\mathbf{x}}
\newcommand{\by}{\mathbf{y}}

\newcommand{\R}{\mathbb R}

\definecolor{NUSBlue}{RGB}{0,61,124}   

\usepackage{tikz}
\usepackage{mathtools}

\def\nudge{.5}

\tikzset{axis/.style={ultra thin, Grey, -latex, shorten <=-\nudge cm, shorten >=-2*\nudge cm}}
\tikzset{line/.style={thick}}

\usepackage{textcomp}
\DeclareMathAlphabet\mathbfcal{OMS}{cmsy}{b}{n}
\usepackage{hyperref}
\usepackage[linesnumbered,vlined,ruled]{algorithm2e}
\usepackage{amsmath,amsfonts,amssymb,amsthm}
\theoremstyle{plain}
\newtheorem{thm}{Theorem}

\newtheoremstyle{cited}%
  {3pt}
  {3pt}
{\itshape}
  {}
  {\bfseries}
  {.}
  {.5em}
  {\thmname{#1} \thmnumber{#2} \thmnote{\normalfont#3}}
\theoremstyle{cited}
\newtheorem{citedthm}[thm]{Theorem}
\newtheorem{citedlem}[thm]{Lemma}
\newtheorem{citedcor}[thm]{Corollary}

\newtheorem{citedprop}[thm]{Proposition}

\usepackage{graphicx,epstopdf,epsf,subfigure} 
\usepackage{cite}
\usepackage{tabularx}
\usepackage{array}
\usepackage{booktabs}
\usepackage{comment}
\usepackage{float}

\begin{document}
\renewcommand*{\Authsep}{, }
\renewcommand*{\Authand}{, }
\renewcommand*{\Authands}{, }
\renewcommand*{\Affilfont}{\normalsize}   
\setlength{\affilsep}{2em}   
\title{Unconstrained optimization using the directional proximal point method}
\date{}
\author[1]{Ming-Yu Chung}
\author[2]{Jinn Ho}
\author[2,*]{Wen-Liang Hwang}
\affil[1]{The department of mathematics, National Taiwan University, Taipei, Taiwan. }
\affil[2]{the Institute of Information Science, Academia 
  Sinica, Taipei 11529, Taiwan.}
\affil[*]{Corresponding author, whwang@iis.sinica.edu.tw}
\maketitle

\begin{abstract}
This paper presents a directional proximal point method (DPPM) to derive the minimum of any $C^1$-smooth function $f$. The proposed method requires a function persistent a local convex segment along the descent direction at any non-critical point (referred to a DLC direction at the point).  The proposed DPPM can determine a DLC direction by solving a two-dimensional quadratic optimization problem, regardless of the dimensionally of the function variables. Along that direction, the DPPM then updates by solving a one-dimensional optimization problem. This gives the DPPM advantage over competitive methods when dealing with large-scale problems, involving a large number of variables. We show that the DPPM converges to critical points of $f$. We also provide conditions under which the entire DPPM sequence converges to a single critical point. For strongly convex quadratic functions, we demonstrate that the rate at which the error sequence converges to zero can be R-superlinear, regardless of the dimension of variables. 
\end{abstract}

\section{Introduction}

The classical approach to unconstrained optimization involves searching for a local minimum without any  knowledge other than that pertaining to function $f$.
In this paper, we make the explicit assumption that $f: \R^n \rightarrow \R \cup \{\infty\} \in \Gamma$, where $\Gamma$ denotes the class of closed, continuously differentiable functions, and $\Omega_0= \{\bx| f(\bx) \leq f(\bx^0)\}$ is a bounded set for a given $\bx^0$. 
First-order methods can be empolyed in the iterative algorithms, wherein the next iterate lies in a direction that capable of decreasing the objective value of the current iterate~\cite{boyd2004convex,bertsekas2014constrained,van2020,pil2020,ruder2016overview}. 
The gradient and proximal point methods are arguably the methods most commonly used for un-constrained optimization.
The gradient method updates iterate $\bx$
using step-size $t > 0$, with the gradient at $\bx$, as follows:
\begin{align*}
\bu = \bx - t \nabla f(\bx).
\end{align*}
The proximal point method (PPM) \cite{trefethen1997numerical,bertsekas2015convex,bauschke2011convex} obtains $\bu$, using the gradient at $\bu$:
\begin{align*}
\bu = \bx -  t \nabla f(\bu).
\end{align*}
The proposed directional proximal point method (DPPM) updates $\bx$ using parameter $t > 0$ and $\bar \bp$ (a descent direction of $f$ at $\bx$\footnote{$\bar \bp$ is a descent direction of $f$ at $\bx$ if and only if $\nabla f(\bx)^\top \bar \bp < 0$.}) to obtain
\begin{align} \label{dppmdefinition}
\bu = \bx - t (\nabla f(\bu)^\top \bar \bp) \bar \bp = \bx + w \bar \bp
\end{align}
where $w \geq 0$ is the step size.

We show that the next iterate of the DPPM is the solution to the following:
\begin{align} \label{dppmwform0}
\arg\min_{w \geq 0}\frac{1}{2t} \|w\|^2 + f(w \bar \bp + \bx) \in \R_+.
\end{align}
This optimization problem resembles the following method, which is used to derive the next iterate of the PPM:
\begin{align} \label{ppmwform0}
\arg\min_{\bu} f(\bu) + \frac{1}{2t} \|\bu - \bx\|^2 \in \R^n.
\end{align}
The main difference between (\ref{dppmwform0}) and (\ref{ppmwform0}) lies in the dimensions of the solution space. The operator of the DPPM (\ref{dppmwform0}) is a one-dimensional optimization problem, involving a search for the optimal step-size $w \in \R_+$. By contrast, the operator of PPM (\ref{ppmwform0}) is a multi-dimensional optimization problem, search for the optimal point in $\R^n$. 

DPPM is applicable if $f \in \Gamma$, and satisfies the assumption of directional local convexity: \\
(DLC) There exists a descent direction $\bar \bp$ at non-critical point $\bx$ of $f$ and   
$v(\bx, \bar \bp) > 0$, such that $f$ is convex over the segment $\{\bu = \bx + w \bar \bp$ with $0 \leq w \leq v(\bx, \bar \bp) \}$. In other words, 
for $w \in [0, v]$, we obtain the following: 
\begin{align} \label{ASdef}
f(\bx + w \bar \bp) \geq f(\bx) + w \nabla f(\bx)^\top \bar \bp.
\end{align}
We define a DLC direction as the direction at any point that satisfies the (DLC). 
The class of smooth functions that holds (DLC) is broad. In fact, it is in the class of $C^1$-smooth, but not strictly concave functions (see Appendix \ref{AS}). 
Various functions satisfying (DCL) are presented in the following:
\begin{itemize}
\item (DLC) holds for any convex function $f \in \Gamma$, where $v(\bx, \bar \bp)$ can be any number in any descent direction $\bar \bp$ of $\bx$.
\item Suppose that Hessian $H$ of a non-convex function $f$ is continuous. Then, (DLC) holds at $\bx$, where $H(\bx)$ has an eigenvector that is not orthonormal to $\nabla f(\bx)$ corresponding to a positive eigenvalue. \\
Without a loss of the generality, we let $\bx = 0$ and $f(0) = 0$. Since the Hessian of $f$ is continuous, for $\bu = w \bar \bp$ with $w > 0$ sufficiently small, we obtain the following:
\begin{align} \label{AS3}
f(w \bar \bp) = w \bar \bp^\top \nabla f(0) + w^2\bar \bp^\top H(0) \bar \bp.
\end{align}
Suppose that $\bar \bp$ is the eigenvector of $H(0)$ corresponding to positive eigenvalue $\lambda^2$.
If $\bar \bp$ and $\nabla f(0)$ are not orthogonal to each other, then  
we can set $\bar \bp$, such that  $\bar \bp^\top \nabla f(0) < 0$. In accordance with (\ref{AS3}),  we obtain the following:
\begin{align*}
f(w \bar \bp) =  w \bar \bp^\top \nabla f(0) + w^2 \lambda^2;
\end{align*}
which means that $f(w \bar \bp)$ is a strictly convex quadratic function along descent direction $\bar \bp$.
\end{itemize}
The contributions of this paper are as follows:
\begin{itemize}
\item We demonstrate that of a given a DCL direction, the optimal step-size (that is, the solution to problem (\ref{dppmwform0})) not only exists but is unique. Regardless of the number of variables in the problem,
we demonstrate that a DLC direction can be derived by solving a quadratic optimization problem with two dual variables and the optimal step-size can be derived by solving the one-dimensional optimization problem. This makes the DPPM competitive for problems involving a large number of variables.

\item The sequence of iterates converges to critical points of $f$, provided that the search direction $\bar \bp_k$ does not tend toward orthogonality to gradient $\nabla f(\bx^k)$. This supposition can be made precise using the following gradient-related definition \cite{Ber08}: $\{\bar \bp_k\}$ is gradient-related to  $\{\bx^k\}$ if  for any subsequence $\{\bx^k\}_{k \in K}$ that converges to a non-critical point of $f$, the corresponding $\{\bar \bp_k\}_{k \in K}$ satisfies 
\begin{align} \label{gradient-related}
\lim\sup_{k \in K, k \rightarrow \infty}\bar \bp_k^\top \nabla f(\bx^k) < 0.
\end{align}
\item We also provide conditions for $\{\bar \bp_k\}$ and the local properties of a critical point of $f$, such that the entire DPPM iterates converges to a single critical point.
\end{itemize}

Studying the rate (speed) of convergence is of practical value as it is often a dominant criterion in algorithm selection for solving a particular problem. In the current paper, we assess the rate of DPPM convergence for convex functions and strongly convex quadratic functions. 
\begin{itemize}
\item For convex functions, the DPPM is as efficient as the PPM in deriving a sub-optimal solution measured in terms of the number of iterations. If fact, the DPPM uses ${\cal O}(\frac{1}{\epsilon})$ iterations to reach an $\epsilon$-suboptimal solution, and this method can be accelerated to  ${\cal O}(\frac{1}{\sqrt\epsilon})$. 
\item For quadratic functions, we show that the rate at which the error sequence of the DPPM converges to zero  is R-linear with rate $(\frac{1}{1+ \lambda m})^\frac{1}{dim}$, where $m$ is the smallest eigenvalue of the quadratic function, $\lambda > 0$, and $dim$ refers to the dimension of the input variables.  Leveraging the fact that $\lambda$ can be any positive number, we demonstrate that it is possible to select sequences of descent directions ($\bar \bp$) and parameters ($t$) for the DPPM, such that the convergence rate of the error sequence to zero is R-superlinear.
\end{itemize}

In contrast to the last contribution, the convergence rate bound for strongly quadratic functions using the steepest descent method with exact line search is Q-linear with rate $\frac{M-m}{M+m}$, where $M$ and $m$  respectively refer to the largest and smallest eigenvalues of a quadratic function with a positive definite Hessian matrix (Theorem 3.3 \cite{nocedal2006numerical}). The Barzilai and Borwei (BB) method \cite{barzilai1988two} can achieve a R-superlinear rate of convergence for two variable cases. Note that in the $n$-dimensional case, the convergence rate of the BB method is R-linear \cite{raydan1993barzilai,dai2002r}. The step-size rule of the BB method was stabilized in \cite{burdakov2019stabilized}, such that long step-sizes causing divergence away from the optimal solution can be avoided.
Note also that given a finite number of iterations, the conjugate gradient method can be used find the optimal solution of a quadratic function. Nonetheless, the success of the conjugate gradient method requires that the search direction and step-size calculations be consistent with data generated by a quadratic function. Any deviation from the quadratic function can seriously degrade performance. As indicated in \cite{fletcher2005barzilai}, this is the situation in which other methods might are worth considering as an alternative to the conjugate gradient method.

The remainder of the paper is organized as follows. Section \ref{backgrounds} reviews related works. Section \ref{sec:DPPM} presents the optimization methods used to derive a DLC direction and the optimal step-size along that direction. Section~\ref{sec:generalDPPM} 
outlines the convergence of the DPPM. 
Section~\ref{convergencerate} analyzes the rate of convergence for convex functions and strongly convex quadratic functions. 
Section~\ref{experiment} outlines the implementations, experiment results for a non-convex function, and the convergence rate of the DPPM for a strongly convex quadratic function. 
Concluding remarks are presented in Section \ref{sec:conclusion}. \\

 {\bf{Notation:}} \\
\noindent$\bullet$ $\|\bx\|$ denotes the 2-norm of $\bx$. \\
 \noindent$\bullet$ $\bar \bp$ is the unit norm directional vector of $\bp$. \\
 \noindent $\bullet$ Boldfaced letters are vectors.

\section{Related works} \label{backgrounds}

Standard first-order methods for unconstrained optimization problems rely on ``iterative descent", which works as follows: Start at point $\bx^0$ and successively generate iterates $\bx^1$, $\bx^2$, $\cdots$, such that $f$  decreases in each iteration. The aim is to decrease $f$ to its minimum.  In many cases, this can be achieved by adopting a strategy accepting a step-size in the search direction according to properties of $f$. 
For a detailed discussion of the topic, see \cite{nocedal2006numerical, Ber08}.
A classical analysis of convergence when using first-order methods for smooth functions reveals the following result: 
\begin{citedprop}(Proposition 1.2.1 \cite{Ber08})
Let $f$ be a smooth function and let $\{\bx^k\}$ be a sequence generated using a gradient-related method $\bx^{k+1} = \bx^k + \alpha_k \bar \bp_k$ where $\{\bar \bp_k\}$ is gradient related (\ref{gradient-related}) and the selection of $\alpha_k$ is based on the minimization rule (exact line search), or Armijo rule. Then, every limit point of $\{\bx^k\}$ is a critical point of $f$.
\end{citedprop}
The exact line search is an effective heuristic used in the selection of the step-size; however, where the number of variables is large, computing the minimization precisely can be too costly for a single step in any gradient-related method.

To avoid the excess computation commonly associated with the minimization rule, it is possible to use inexact line searches based on successive step-size reduction. 
The Armijo rule \cite{armijo1966minimization} establishes convergence to critical points of smooth functions using an inexact line search with a simple ``sufficient decrease" condition. This condition ensures that the line search step is not too large. Armijo back-tracking begins with a fixed initial step-size and geometrically decreases the step-size until the Armijo condition is established. Note that this is perhaps the most common line search method in practice. The approximately exact line search algorithm \cite{fridovich2020approximately} uses solely function evaluations to minimize $f$, which is assumed to be strongly convex with a Lipschitz continuous gradient. 
The algorithm iteratively decreases and/or increases the step-size along a descent direction in order to identify a step-size within a constant fraction of the exact line search minimizer (without exceeding it). Note that the step-size does not necessarily satisfy the Armijo condition. 

The first-order methods have recently been challenged, as the number of variables $n$ in many optimization problems has grown far too large to apply methods that require more than ${\cal O}(n)$ operations per iteration \cite{asl2020analysis,fletcher2005barzilai}. If efficiency is the sole concern for a large-scale optimization problem, then the BB method \cite{barzilai1988two} is an option, due to the fact that it requires only ${\cal O}(n)$ floating point operations and a gradient evaluation for an update. In fact, the BB method is globally convergent at a R-linear convergent rate for strongly quadratic functions of any dimension. The search direction of the method is always along the direction of the negative gradient. The step-size is not the conventional choice for the steepest descent method and it does not use back-tracking line search and therefore cannot guarantee a descent in each update. Note however that enforcing descent in the negative gradient direction in every iteration can destroy some of the local properties of a function, such that this method becomes just another version of the (slow) steepest descent method. 
Raydan \cite{raydan1997barzilai} established the global convergence of the BB method for non-quadratic functions by incorporating the non-monotone line search of Grippo, Lampariello, and Lucidi \cite{grippo1986nonmonotone}. Since that time, this method has been extended to many fields \cite{zhou2006gradient}, such as convex constrained optimization, non-linear systems, and support vector machines.

\section{Characterizing DPP-updates via optimization}  \label{sec:DPPM}

The DPP-update involves two tasks: Determining the direction satisfying the (DLC) and determining the optimal step-size along a DLC direction. We show that both tasks can be accomplished by solving low-dimensional optimization problems.

\subsection{Optimal step-size}

DPPM updates iterate $\bx$ when $t >0$ using 
\begin{align} \label{updatedppm}
\bu = \bx - [t \bar \bp^\top \nabla f(\bu)] \bar \bp = \bx + w \bar \bp,
\end{align}
where $\bar \bp$ is a direction of descent in which $f$ exhibits local convex behavior at $\bx$ and $w \geq 0$ is the step-size.
The monotonicity of the convex function for any $\bu_1$ and $\bu_2$
\begin{align} \label{monoprop}
[\nabla f(\bu_1) - \nabla f(\bu_2)]^\top (\bu_1 - \bu_2) \geq  0,
\end{align} 
holds locally along the ray at $\bar \bp$.
\begin{citedlem}\label{dppmmonoprop}
Suppose that $f \in \Gamma$ satisfies (DLC) . Let $\bu_1 = \bx + w_1 \bar \bp$ and $\bu_2 = \bx + w_2 \bar \bp$ for $w_1, w_2 \in [0, v(\bx, \bar \bp)]$ (the convex segment of $f$ at $\bx$ along $\bar \bp$). Then,
\begin{align} \label{dppmmonoprop1}
[\nabla f(\bu_1)^\top \bar \bp - \nabla f(\bu_2)^\top \bar \bp]\bar \bp^\top (\bu_1 - \bu_2) \geq 0.
\end{align}
\end{citedlem}
\proof
See Appendix \ref{lemma-1}
\qed

\begin{citedcor}\label{dppmmonocor}
$\bar \bp^\top \nabla f(\bx + w \bar \bp) \leq 0$ is an increasing function of $w \in [0, v(\bx, \bar \bp)]$ ($|\bar \bp^\top \nabla f(\bx + w \bar \bp)|$ is thus a decreasing function of $w$). Therefore, for $w \in [0, v(\bx, \bar \bp)]$, 
\begin{align} \label{monodir}
| \bar \bp^\top \nabla f(\bx) | \geq |\bar \bp^\top \nabla f(\bx + w \bar \bp)|
\end{align}
and $w + t \bar \bp^\top \nabla f(\bx + w \bar \bp)$ is a strictly increasing function. 
\end{citedcor}

The following property for function $f$ consisting of local strictly convex behavior can be derived using the strict monotonicity property with derivation in parallel with Lemma \ref{dppmmonoprop}.

\begin{citedcor} \label{sdppmmonoprop}
Suppose that the assumption pertaining to Lemma \ref{dppmmonoprop} holds and that the function over ray $\bar \bp$ at $\bx$ is locally strictly convex. Then, 
\begin{align} \label{sdppmmonoprop1}
[\nabla f(\bu_1)^\top \bar \bp - \nabla f(\bu_2)^\top \bar \bp]\bar \bp^\top (\bu_1 - \bu_2) > 0.
\end{align}
Hence, $\bar \bp^\top \nabla f(\bx + w \bar \bp) < 0$ is a strictly increasing function of $w \in [0, v(\bx, \bar \bp)]$ ($|\bar \bp^\top \nabla f(\bx + w \bar \bp)|$ is thus strictly decreasing function of $w$). Therefore, for $w \in [0, v(\bx, \bar \bp)]$, 
\begin{align} \label{smonodir}
| \bar \bp^\top \nabla f(\bx) | > |\bar \bp^\top \nabla f(\bx + w \bar \bp)|.
\end{align}
\end{citedcor}

We consider the Moreau envelope function~\footnote{$\min_{\bu} \frac{1}{2t} \| \bu - \bx \|^2 + f(\bu)$} with constraints $\bu = \bx + w \bar \bp$ and $w \geq 0$:
\begin{align} \label{dppmoptfor}
\begin{cases}
\min_{w,\bu} \frac{1}{2t} \| \bu - \bx\|^2 + f(\bu) \\
\bu - \bx = w\bar \bp \\
w \geq 0.
\end{cases}
\end{align}
Substituting $\bu - \bx = w\bar \bp$ into the objective function yields the following envelope function for the DPPM:
\begin{align} \label{dppmwform}
\min_{w \geq 0}\frac{1}{2t} \|w\|^2 + f(w \bar \bp + \bx).
\end{align}
Taking the derivative of (\ref{dppmwform}) with respect to $w > 0$  and setting the result to zero, we obtain the following:
\begin{align} \label{dppmwvalue}
\frac{1}{t} w^* + \bar \bp^\top \nabla f(\bx + w^* \bar \bp) = 0. 
\end{align}

Below, we show that the DPP-update (\ref{updatedppm}) can be obtained by solving the optimization problem presented in (\ref{dppmwform}).

\begin{citedlem}\label{dppmunique}
Suppose that $f \in \Gamma$ holds (DLC) at $\bx$ with $w \in [0, v(\bx, \bar \bp)]$, where $\bar \bp$ is in the  descent direction of $f$ at $\bx$. If we choose $t \leq \frac{v(\bx, \bar \bp)}{| \bar p^\top \nabla f(\bx)|}$ and let  $w^*$ denote a solution to (\ref{dppmwform}), then  \\
(i)  $w^* \in (0, v(\bx, \bar \bp)]$ where
\begin{align} \label{dppmupdate0}
\bu = \bx  + w^* \bar \bp = \bx - t \bar \bp \bar \bp^\top \nabla f(\bu),
\end{align}
(ii) the solution $w^*$ to (\ref{dppmwform}) is unique. \\
 \end{citedlem}
\proof

(i) In accordance with variational inequality \cite{kinderlehrer2000introduction}, $w^*$ is a solution to (\ref{dppmwform}) if and only if
\begin{align} \label{dppmvarI}
(\frac{1}{t} w^* + \bar \bp^\top \nabla f(\bu)) (w - w^*) \geq 0 \text{ for all $w \geq 0$}.
\end{align}
This equation is equivalent to    
\begin{align} \label{dppmvi}
\begin{cases}
\frac{1}{t} w^* + \bar \bp^\top \nabla f(\bu) = 0 \text{ if $w^* > 0$} \\
 \bar \bp^\top \nabla f(\bu) \geq 0 \text{ if $w^* = 0$}.
\end{cases}
\end{align}
From (\ref{dppmvi}) and $\bu = \bx + w^* \bar \bp$, we can assert that
\begin{align} \label{dppminnerproduct1}
\begin{cases}
\bar \bp^\top \nabla f(\bu) < 0 \text{ if $w^* > 0$} \\
\bar \bp^\top \nabla f(\bx) = 0  \text{ if $w^* = 0$}.
\end{cases}
\end{align}
The second part of (\ref{dppminnerproduct1}) violates the assumption that $\bar \bp$ is in the descent direction at $\bx$. 
From (\ref{dppmwvalue}), (\ref{dppminnerproduct1}), Corollary \ref{dppmmonocor}, and $t \leq \frac{v(\bx, \bar \bp)}{| \bar \bp^\top \nabla f(\bx)|}$, we obtain the following:
\begin{align}
w^*  & =  - t \bar \bp^\top \nabla f(\bx + w^* \bar \bp)  \nonumber  \\
  & \leq t | \bar \bp^\top \nabla f(\bx) |  \leq v(\bx, \bar \bp) \label{tparameter}
\end{align}

(ii) We prove that the solution is unique via contradiction. Suppose that $w_1$ and $w_2 \neq w_1$ are solutions with 
 \begin{align*}
\bu_1 = \bx - t \bar \bp \bar \bp^\top \nabla f(\bu_1) = \bx + w_1 \bar \bp
\text{ and }
\bu_2 = \bx - t \bar p \bar \bp^\top \nabla f(\bu_2) = \bx + w_2 \bar \bp.
\end{align*}
Simple algebra yields
\begin{align*}
\bu_1 - \bu_2 = t \bar \bp [\bar \bp^\top (\nabla f(\bu_2) - \nabla f(\bu_1))] = (w_1 - w_2) \bar \bp.
\end{align*}
From the above and Lemma \ref{dppmmonoprop}, we obtain 
\begin{align*}
 [\bar \bp^\top (\nabla f(\bu_1) - \nabla f(\bu_2))]\bar \bp^\top (\bu_1 - \bu_2) = - t \|\bar \bp^\top (\nabla f(\bu_2) - \nabla f(\bu_1))\|^2 = \frac{-(w_1 - w_2)^2}{t} \geq 0,
\end{align*}
which contradicts the assumption that $w_1 \neq w_2$.

\qed

\begin{citedcor} \label{convexpara}
Suppose that $f \in \Gamma$ and $f$ is a convex function. If $| \bar \bp^\top \nabla f(\bx)|$ is bounded for all $\bx$ , then $t$ can be a constant. 
\end{citedcor}

\subsection{Direction of DLC (DLC direction)} \label{sec:DLCdirection}

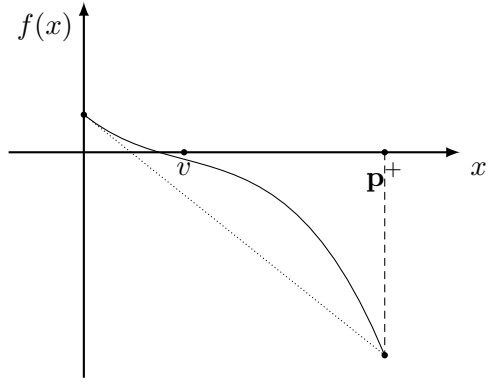
\begin{figure} 
\centering
\begin{tikzpicture}
\draw[-latex,thick] (-1,0) -- (5,0) node[below right] {$x$};
\draw[-latex,thick] (0,-3) -- (0,2) node[below left] {$f(x)$};

\draw[densely dotted](0, 0.5) -- (4, -2.7);
\draw[densely dashed](4,0) -- (4,-2.7);

\fill (0, 0.5) circle (1.2 pt);
\fill (4, -2.7) circle (1.2 pt);

\fill (4/3,0) circle (1.2 pt) node[below]{$v$};
\fill (4,0) circle (1.2 pt) node[below]{$\mathbf{p}^+$};

\draw[domain= 0:4, yscale = 0.1] plot (\x,{(4*(\x-1)^2 -\x^3 +1 )});
\end{tikzpicture}
\caption{The direction satisfying the first constraint in (\ref{DLCdirection}) implies that there exists an convex segment ($[0,v]$) along the direction, due to the fact that $f$ is not a strictly concave function at $\bx$. } \label{DLCfig}
\end{figure}

The DLC direction $\bar \bp$ at $\bx$ is a direction satisfying (DLC); hence, it is a descent direction possessing a convex segment $[\bx,  \bx+ v \bar \bp]$ with $v > 0$ (that is, for any $\bu, \bv \in [\bx,  \bx+ v \bar \bp]$, $f(\bu) \geq f(\bv) + \nabla f(\bv)^\top (\bu-\bv)$). 
In accordance with Corollary \ref{characterizeDLC}, $f$ satisfying (DLC) is not strictly concave at $\bx$. For the purpose of seeking the DLC direction $\bar \bp^*$ at $\bx$, we can formulate the following optimization problem:
\begin{align} \label{DLCdirection}
\begin{cases}
\min_{\bp}\frac{1}{2} \| \bp \|^2  \\
 \text{ subject to }   f(\bx + \bp) \leq f(\bx) + \nabla f(\bx)^\top \bp \text{ and }  \nabla f(\bx)^\top \bp \leq 0. 
\end{cases}
\end{align}
The first constraint and the objective function indicate that any $\bp = r  \bar \bp^*$ with $\br \leq \|\bp^*\|$ will not satisfy the inequality. 
Therefore, if $\bar \bp^* \neq 0$, then $\bx+ [0,  \| \bp^*\|] \bar\bp^*$ contains a convex segment along direction $\bar \bp^*$ at $\bx$. As shown in Figure \ref{DLCfig} is the direction satisfying the constraint which possess a convex segment $[0, v]$. The second constraint indicates that $\bar \bp^*$ is a descent direction. Since  (\ref{DLCdirection}) has a trivial solution at $\bp = 0$ (indicating that the assumption of DLC is not satisfied at $\bx$), the first constraint of (\ref{DLCdirection}) is amended to remove $\bp = 0$ as a solution by introducing $\delta > 0$ with 
\begin{align}
f(\bx + \bp) + \delta \leq f(\bx) + \nabla f(\bx)^\top \bp.
\end{align}
Note that $\delta$ is set a a small value.

If we let $g_0(\bp)$, $g_1(\bp)$, and $g_2(\bp)$ at $\bx$ respectively be
\begin{align} 
g_0(\bp) & = \frac{1}{2} \|\bp\|^2, \label{auxiliaryfun1}\\
g_1(\bp) & =  f(\bx + \bp) + \delta - f(\bx) - \nabla f(\bx)^\top \bp \leq 0, \label{auxiliaryfun2}\\
g_2(\bp) & =   \nabla f(\bx)^\top \bp \leq 0, \label{auxiliaryfun3}
\end{align}
then (\ref{DLCdirection}) becomes
\begin{align} \label{DLCdirection1}
\begin{cases}
\min g_0(\bp)   \\
 \text{ subject to }   g_1(\bp) \leq 0  \text{ and }  g_2(\bp) \leq 0.
\end{cases}
\end{align}
We propose solving this problem using penalty with the linearization method \cite{pshenichnyj2012linearization}, which transforms (\ref{DLCdirection1}) into the following quadratic programming problem, where $\bp_0$ is the current estimate:
\begin{align} \label{DLCdirection2}
\begin{cases}
\min_{\bp}\nabla g_0(\bp_0)^\top (\bp- \bp_0) + \frac{\mu}{2} \| \bp -\bp_0 \|_2^2   \\
 \text{ subject to }   g_i(\bp_0) +\nabla g_i(\bp_0)^\top (\bp - \bp_0) \leq 0  \text{ for $i = 1, 2$ }.
\end{cases}
\end{align}
The penalty $\frac{\mu}{2} \|\bp- \bp_0\|^2$ in the objective is meant to penalize large deviations from $\bp_0$. 
The Lagrangian function of  (\ref{DLCdirection2}) is 
\begin{align} \label{Lagrangian}
L_{\mu}(\bp, \lambda_i) = 
\nabla g_0(\bp_0)^\top (\bp- \bp_0) + \frac{\mu}{2} \| \bp -\bp_0 \|_2^2  + \sum_{i=1,2} \lambda_i [g_i(\bp_0) +\nabla g_i(\bp_0)^\top (\bp - \bp_0)]
\end{align}
where $\lambda_i \geq 0$. The KKT condition for convex analysis indicates that $\bp$ is a solution to (\ref{DLCdirection2}) if and only if there exists $\lambda_i \geq 0$ such that
\begin{align} \label{KKT1}
\begin{cases}
\nabla_{\bp} L_{\mu}(\bp, \lambda_i) = \nabla g_0(\bp_0) + \mu(\bp-\bp_0) + \sum_i \lambda_i \nabla g_i(\bp_0) = 0 \\
\lambda_i (g_i(\bp_0) +\nabla g_i(\bp_0)^\top (\bp - \bp_0)) = 0.
\end{cases}
\end{align}

The following lemma indicates that we can solve (\ref{DLCdirection2}) to obtain a solution that satisfies the KKT condition for (\ref{DLCdirection1}).
\begin{citedlem} \label{DLClemma}
$\bp_0$ satisfies the necessary condition for the minimization of (\ref{DLCdirection1}) if and only if $\bp =\bp_0$ is a solution to (\ref{DLCdirection2}).
\end{citedlem}
\proof
Suppose that  $\bp_0$ satisfies the necessary condition for the minimization of (\ref{DLCdirection1}). Then there must exist $\lambda_i \geq 0$, such that
\begin{align}
\nabla g_0(\bp_0) + \sum_{i=1,2} \lambda_i \nabla g_i(\bp_0) = 0 \text{ and } \lambda_i g_i(\bp_0) = 0.
\end{align}
Comparison with (\ref{KKT1}), $\bp = \bp_0$ is also the solution of (\ref{DLCdirection2}).

If we suppose that $\bp = \bp_0$ is a solution to (\ref{DLCdirection2}), then from (\ref{KKT1}), we obtain 
\begin{align} \label{KKT0}
\begin{cases}
\nabla g_0(\bp_0)  + \sum_i \lambda_i \nabla g_i(\bp_0) = 0 \\
\lambda_i g_i(\bp_0) = 0.
\end{cases}
\end{align}
This is the necessary condition for the minimization of (\ref{DLCdirection1}).

\qed

It is preferable that the solution to (\ref{DLCdirection1}) be derived from the dual problem of (\ref{DLCdirection1}), due to the fact that the dual function consists of only two variables, regardless of the dimensionality of the primal variables.
The dual problem of (\ref{DLCdirection2}) is calculated as follows: 
\begin{align} \label{dual0}
\max_{\lambda_i \geq 0} \min_{\bp} L_{\mu}(\bp, \lambda_i) = \max_{\lambda_i \geq 0} d(\lambda_i),
\end{align}
where the Lagrangian function $L_{\mu}$ is given in (\ref{Lagrangian}) and $d$ is the dual function. From (\ref{KKT1}), we obtain
\begin{align} \label{sol}
\bp = \bp_0 -  \frac{1}{\mu}( \nabla g_0(\bp_0) +\sum_i \lambda_i \nabla g_i(\bp_0)).
\end{align}
Substituting this into function $L_{\mu}$, we obtain the following two-variable quadratic function:

\begin{align}
d(\lambda_i) 
&=
[\sum_{i = 1,2} \lambda_i g_i(\bp_0)] \nonumber
- 
\frac{1}{2\mu} \|\nabla g_0(\bp_0) + \sum_{i = 1,2} \lambda_i g_i(\bp_0)\|^2
\label{prop_Llam}
\end{align}

Let $\lambda_i^+ \geq 0$ with $i=1,2$ be the Lagrangian multipliers used to optimize the dual function with simple constraint:
\begin{align*}
\max_{\lambda_i \geq 0} d(\lambda_i). 
\end{align*}
By substituting $\lambda_i^+$ into (\ref{sol}), the solution to (\ref{DLCdirection2}) is 
\begin{align}
\bp^+ = \bp_0 -  \frac{1}{\mu}( \nabla g_0(\bp_0) + \sum_i \lambda_i^+ \nabla g_i(\bp_0)).
\end{align}
The above procedure can be used repeatedly to solve (\ref{DLCdirection2}) where $\bp_0 \leftarrow \bp^+$ ($\bp^+$ is the update of $\bp_0$) until $\bp^+$ is sufficiently close to $\bp_0$, thereby providing an approximate solution to (\ref{DLCdirection1}) in accordance with Lemma \ref{DLClemma}.

\section{Convergence analysis}\label{sec:generalDPPM}

We posit that if search direction $\bar \bp_k$ does not become orthogonal to gradient $\nabla f(\bx^k)$, then the DPPM converges to critical points of $f$. We also present the conditions under which the DPPM converges to a single critical point of $f $. 

\subsection{Convergence to critical points}

Conventional analysis has demonstrated that any smooth, bounded below function possessing step-sizes satisfying the Armijo condition \cite{nocedal2006numerical}; that is, there exists $0 < c_1 < 1$, such that
\begin{align} \label{dppmwolfe}
f(\bx^{i+1}) \leq f(\bx^i) + c_1 w_k \nabla f(\bx^i)^\top \bar \bp_i \; \;\forall i,
\end{align}
and the convergence of gradient-related methods in the selection of step-sizes satisfying the Armijo condition \cite{armijo1966minimization}.
The Armijo condition is also referred to as the sufficient decrease condition, which ensures a sufficient decrease of $f$ in each update, based on the fact that an insufficient reduction in $f$ in each step could cause an algorithm to fail in its convergence to the minimizer of $f$.

\begin{citedlem} \label{dppmlem}
Suppose that $f \in \Gamma$ holds (DLC). Let 
\begin{align*}
\bx^{k+1} = \bx^k - t_k (\nabla f(\bx^{k+1})^\top \bar \bp_k) \bar \bp_k = \bx^k + w_k \bar \bp_k,
\end{align*}
and $w_k \in [0, v(\bx^k, \bar \bp_k)]$ for all $k$. 
The DPPM is a descent algorithm with 
\begin{align} 
f(\bx^{k+1}) & \leq f(\bx^{k}) - t_k | \nabla f(\bx^{k+1})^\top \bar \bp_k|^2  \label{dppmdescent}\\
& =  f(\bx^k) + w_k \nabla f(\bx^{k+1})^\top \bar \bp_k.  \label{dppmwolfe1}
\end{align}
Note that (\ref{dppmwolfe1}) is obtained using $w_k = t_k|\nabla f(\bx^{k+1})^\top \bar p_k|$ (\ref{dppmwvalue}) and $f(\bx^{k+1})^\top \bar \bp_k < 0$.
\end{citedlem}
\proof
See Appendix \ref{DPPMdescent}.
\qed

Due to the fact that $| \bar \bp_k^\top \nabla f(\bx^k) | > | \bar \bp_k^\top \nabla f(\bx^{k+1}) |$ (\ref{monodir}), we can introduce $0 < c_1(i) \leq 1$ in order to express (\ref{dppmwolfe1}) as
\begin{align} \label{dppmarmijo}
f(\bx^{i+1}) \leq f(\bx^i) + c_1(i) w_i \nabla f(\bx^i)^\top \bar \bp_i.
\end{align}
Clearly, this inequality remains true if $c_i(i)$ is replaced with any scalar in $(0,c_1(i)]$.  This allows us to define the Amijo parameter up to $k$ as $c_{1,k}$ where
\begin{align*}
\cap_{i \leq k} (0, c_1(i)] = (0, c_{1, k}]
\end{align*}
and obtain 
\begin{align} \label{Amijok}
f(\bx^{i+1}) \leq f(\bx^i) + w_i c_{1, k} \bar \bp_i^\top \nabla f(\bx^i) \text{ for $i \leq k$}.
\end{align}
Note that (\ref{Amijok}) indicates that the Amijo condition (\ref{dppmwolfe}) holds for all $i \leq k$ with parameter $c_{1,k}$. Extending $k \rightarrow \infty$ requires the following technical lemma.

\begin{citedlem}\label{sdppmlem}
Suppose that there exists $k$, such that $f$ over ray $\bar \bp_k$ at $\bx^k$ presents strictly local behavior over $(0, v(\bx^k, \bar \bp_k)]$. Then, $\lim_{k\rightarrow \infty} c_{1, k} = c_1 \in [0, 1)$. 
\end{citedlem}
\proof

Following $| \bar \bp_k^\top \nabla f(\bx^k) | > | \bar \bp_k^\top \nabla f(\bx^{k+1}) |$ (\ref{smonodir}),  there exists $0 < c_1(k) < 1$, such that
\begin{align} \label{dppmwolfe2}
c_1(k) | \bar \bp_k^\top \nabla f(\bx^k) | = | \bar \bp_k^\top \nabla f(\bx^{k+1}) |.
\end{align}
Substituting (\ref{dppmwolfe2}) into (\ref{dppmwolfe1}) and using the fact that $\bar \bp_k$ is a descent direction at $\bx^k$, we obtain 
\begin{align*}
f(\bx^{k+1}) \leq f(\bx^k) + w_k c_1(k) \bar \bp_k^\top \nabla f(\bx^k).
\end{align*}
This inequality clearly holds for any in $(0,c_1(k)]$. The Amijo parameter up to $k$ is
\begin{align*}
\cap_{i \leq k} (0, c_1(i)] = (0, c_{1, k}] \subseteq (0, 1).
\end{align*}
The fact that $0 < c_{1, k} < 1$ and $c_{1, k} \leq c_{1, l}$ with $k \geq l$ implies that $\{c_{1,k}\}_k$ is a bounded decreasing sequence with $\lim_{k\rightarrow \infty} c_{1, k} = c_1 \in [0, 1)$.

\qed

The fact that $\{\bar \bp_k\}$ is gradient-related to  $\{\bx^k\}$ ensures that  
if sequence $\{\nabla f(\bx^k)\}_{k \in K}$ tends toward a non-zero vector, the corresponding subsequence of $\bar \bp_k$ will not tend to be orthogonal to $\nabla f(\bx^k)$\cite{nocedal2006numerical,Ber08}. Gradient-related directions include the negative gradient direction or $-D_k \nabla f(\bx^k)$, where $D_k$ is a positive definite matrix bounded above and away from zero \cite{Ber08}. 
An important consequence of convergence is the fact that if $\{\bar \bp_k\}$ is gradient-related and if we use the minimization rule, the limited minimization rule, the Armijo rule, the Goldstein rule, then all limit points of $\{\bx^k\}$ are critical points \cite{Ber08}. The following shows that under mild conditions, convergent is retained for DPPM.

\begin{figure}
\centering
\begin{tikzpicture}
\draw[-latex,thick] (-1,0) -- (6.5,0) node[below right] {$\alpha$};
\draw[-latex,thick] (0,-1) -- (0,4.5) node[below left] {$\phi(\alpha)$};

\draw[densely dotted](0,2.8049) -- (4.2,0.5090);
\draw[densely dashed](4.2,0.5090) -- (4.2,0);

\fill (0,2.8049) circle (1.2 pt);
\fill (4.2,0.5090) circle (1.2 pt);

\fill (0.5,0) circle (1.2 pt) node[below]{$w$};
\fill (2.8,0) circle (1.2 pt) node[below]{$v$};
\fill (4.8,0) circle (1.2 pt) node[below]{$\hat{w}$};
\fill (4.2,0) circle (1.2 pt) 
node[below]{$\alpha'$};

\draw[domain=0:6] plot (\x,{0.5+1/4*sin((\x-4) r)*exp{(2.5-\x)}});
\end{tikzpicture}
\caption{$f(\bx) + \alpha c_1 \nabla f(\bx)^\top \bar \bp$ (dashed line) intersects $\phi(\alpha)$ at $\alpha'$. Note that $\phi(\alpha)$ over $[0, v]$ is convex and $w$ is the step-size to update $\bx$.} \label{mvt}
\end{figure}
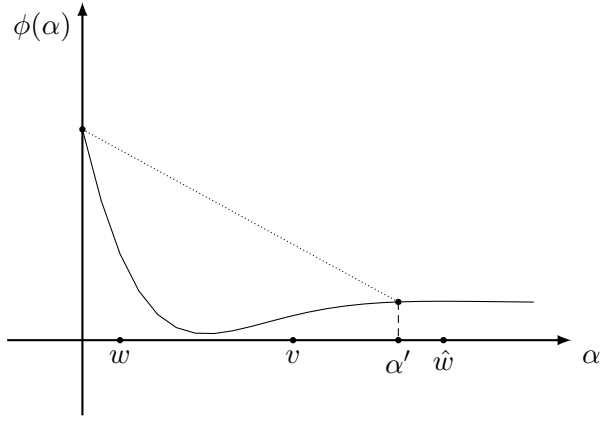

The fact that $\phi(\alpha_k) = f(\bx^k + \alpha_k \bar \bp_k)$ is bounded below for all $\alpha_k > 0$ and line $f(\bx^k) + \alpha_k c_1 \nabla f(\bx^k)^\top \bar \bp_k$ is unbounded below implies that the line must intersect the graph of $\phi(\alpha_k)$ at least once (see Figure \ref{mvt}). Let $\alpha_k' > 0$ be the smallest intersecting value of $\alpha_k$ where
\begin{align} \label{alpha1}
f(\bx^k + \alpha_k' \bar \bp_k) = f(\bx^k) + \alpha_k' c_1 \nabla f(\bx^k)^\top \bar \bp_k.
\end{align}
Let $\{(\alpha_k', w_k)\}$ be a sequence with $\alpha_k' \geq w_k$. 
It is clear there exist $\eta_k \in (0, 1)$ and positive integer $l(k) \geq 1$, such that $\frac{w_k}{l(k) \eta_k} \geq \alpha_k'$. To obtain convergence, we need the following assumption on $\eta_k$: \\
\noindent (AS1) there exists $\eta$ such that $\displaystyle \liminf_{k \rightarrow \infty}\eta_k \geq \eta > 0$.
\begin{citedthm} \label{dppmlimitpoints}
Suppose that $f \in \Gamma$ holds (DLC). Let $\{ \bx^k \}$ and $\{ w_k \}$be the sequences of iterates and step-sizes generated by the DPPM and $\{\bp_k\}$ be gradient-related (\ref{gradient-related}). 
We further suppose that the assumptions pertaining to Lemma \ref{sdppmlem} and (AS1) hold true. Then, \\
(i) Whose sequence $\{ f(\bx^k)\}$ converges to a finite value of $f$. \\
(ii) Every limit point of $\{ \bx^k \}$ is a critical point of $f$. 
\end{citedthm}

\proof

See Appendix \ref{dppmlimitpointsproof}.

\qed

\subsection{Convergence to a single critical point}

A line search algorithm on non-convex functions cannot guarantee global convergence; that is, convergence of the entire sequence to a single critical point from any initial point. 
In \cite{bolte2014proximal}, the global convergence of a sequence for optimization is achieved based on the KL-property of functions. The following demonstrates a different condition under which the entire sequence of the DPPM converges to a single critical point.
Suppose that $\bar \bp$, $\overline{\nabla f(\bu)}$, and $\bar \bp^*$ (that is, $\overline{\bx^* - \bx})$ respectively denote the descent direction at iterate $\bx$, the gradient direction at next iterate $\bu$, and the direction at $\bx$ pointing to a critical point $\bx^*$ of $f$. Let $\bar \bv = \mathbb P(\bar \bp^*)$ denote the orthogonal projection of $\bar \bp^*$ to the plane spanned by $\bar \bp$ and $\overline{\nabla f(\bu)}$. Then, we obtain the following:
\begin{align} \label{gradientdecomp}
\overline{\nabla f(\bu)}= \alpha \bar \bp + \beta  \bv =  \alpha \bar \bp + \beta \bar \bv \| \bv\|.
\end{align}
If we apply inner product $\bar \bp$ to both sides of this equation, we obtain 
\begin{align*}
0 \geq \bar \bp^\top \overline{\nabla f(\bu)} = \alpha + \beta (\bar \bp^\top \bar \bv) \| \bv\|= \alpha + \beta \cos(\bar \bp, \bar \bv) \| \bv\|.
\end{align*}
The above inequality is in accordance with Corollary \ref{dppmmonocor}. If we choose a descent direction, such that $\cos(\bar \bp, \bar \bv) > 0$ and $\overline{\nabla f(\bu)}^\top \bar \bp > \alpha $, then we obtain
\begin{align} \label{cosine}
0 \leq \beta \leq  \frac{-\alpha}{cos(\bar \bp, \bar \bv) \| \bv\|}.
\end{align}

We let $C(\bx^0)$ denote the set of critical points of $f$ to which the sequence of DPPM iterates $\{\bx^k\}$ converges, when the initial point is $\bx^0$. Recall that the sequence of iterates eventually satisfies Fejer monotone with respect to $C(\bx^0)$ \cite{bauschke2011convex} is defined as there exists $k_0$, such that
\begin{align} \label{fejer0}
\| \bx^{k+1} - \bx^* \| \leq \|\bx^k - \bx^*\| \text{ for all $\bx^* \in C(\bx^0)$ and $k \geq k_0$}.
\end{align}
For a sufficiently large $k$, each iterate in the sequence is not strictly farther than its predecessor from any point in $C(\bx^0)$, which means that the norm sequence $\{\| \bx^{k} - \bx^* \|\}$ converges for all $\bx^* \in C(\bx^0)$.

\begin{citedcor} \label{dppmfejer}
Suppose that the assumption pertaining to Theorem \ref{dppmlimitpoints} holds. Let $C(\bx^0)$ denote the set of critical points to which iterates $\{\bx^k\}$ with initial $\bx^0$ converge. $C(\bx^0)$ is not an empty set with $\bx^* \in C(\bx^0)$. Suppose that there is a convex region $\mathcal D$ in the neighborhood $\mathcal N$ of $\bx^*$, such that $\bx^* \in \mathcal D \subseteq \mathcal N$ and $f$ is convex over $\mathcal D$. We further suppose that the descent directions of the DPPM satisfy (\ref{cosine}) and $\bx^k \in \mathcal D$ for $k \geq k_0$.
Then, \\
(i) $\{\bx^k\}$ is a Fejer monotone with respect to $C(\bx^0)$ for $k \geq k_0$. \\
(ii) The entire sequence $\{\bx^k\}$ generated by the DPPM converges to a single critical point of $f$.
\end{citedcor}

\proof
(i) See Appendix \ref{dppmfejerproof}.

(ii) This can be found in \cite{bauschke2011convex}, and we include it here for the sake of completeness. Suppose that sub-sequences $\{\bx^{y(k)}\}_k$ and $\{\bx^{z(k)}\}_k$ respectively converge to critical points $\by$ and $\bz$. In accordance with (i), sequences $\{\|\bx^k - \by\| \}$ and  $\{\|\bx^k - \bz\| \}$ converge. From 
\begin{align*} 
2 (\bx^k)^\top (\bz-\by) = \|\bx^k - \by \|^2 - \| \bx^k - \bz\|^2  + \|\bz\|^2 - \|\by\|^2 \text{ for all $k$,}
\end{align*}
we can deduce  that $\{(\bx^k)^\top (\bz-\by)\}$ also converges, say $\{(\bx^k)^\top (\bz-\by)\} \rightarrow l$. Proceeding to the limits along $\{\bx^{y(n)}\}$ and $\{\bx^{z(n)}\}$ respectively yields
\begin{align} \label{entire1}
l = \by^\top (\bz - \by)
\end{align}
and 
\begin{align} \label{entire2}
l = \bz^\top (\bz - \by).
\end{align}
From (\ref{entire1}) and (\ref{entire2}), $\| \by - \bz \|^2 = 0$; therefore, $\by = \bz$. 

\qed 

The critical points of a non-convex function are usually unions of connected, compact regions.
This corollary indicates that if a DPPM sequence approaches critical points from a convex sub-neighborhood of points and the function restricted to that region is convex, then the entire sequence converges to a single critical point.

\section{Rate of convergence} \label{convergencerate}

Generally, the rate at which the DPPM converge cannot be derived without imposing assumption beyond (DLC), such as the Hessian function of $f$. In this section, we examine the rate at which the DPPM converges for convex functions and strongly convex quadratic functions, which can be used in estimation of local convergence rates for smooth functions.
Analysis of the local rate of convergence describes the local behavior of a given method in the neighborhood of a critical point; however, it disregards the behavior of that method at a distance from the critical point. 

\subsection{Convex functions}

A function satisfying (DLC) assumes that the function exhibits local convex behavior along a descent direction at any non-critical point. 
Obviously, convex function $f \in \Gamma$ must satisfy (DLC) at any non-critical point $\bx$ along any descent direction $\bar \bp$ at any interval $[\bx, \bx + v(\bx, \bar \bp) \bar \bp]$ (i.e., $v(\bx, \bar \bp)$ can be any non-negative real number). This fact is touched on in Lemma \ref{convexpara} which justifies taking constant $t_k$ for all $k$. 
The following reveals that the DPPM is as efficient as the PPM in reaching a sub-optimal solution for a convex function.
\begin{citedlem} \label{dppmexactconvex}
Suppose that $f \in \Gamma$ is convex with an optimal value attained at $f^*$. Let $\{\bar \bp_k\}$ be gradient-related and $t_k = \lambda$ for all $k$. The DPPM achieves an $\epsilon$-suboptimal solution (i.e., $f(\bx^k) -f^* \leq \epsilon$) using $k$ of iterations in the order of ${\cal O}(\frac{1}{\epsilon})$. 
Furthermore, the number of iterations can be reduced to ${\cal O}({\frac{1}{\sqrt{\epsilon}}})$.
\end{citedlem}
\proof
our analysis is similar to that for the PPM \cite{van2020}. 
Using $t_k = \lambda$, we obtain the following:
\begin{align} 
2k(f(\bx^k) - f^*) \lambda & \leq 
2k(f(\bx^k) - f^*) \lambda + \sum_{i=0}^{k-1} w_i^2 \\ 
&\leq \sum _{i=0}^{k-1} 2\lambda  (f(\bx^{i+1}) - f(\bx^*)) + \sum_{i=0}^{k-1} w_i^2  \nonumber \\
&  \leq \sum_{i=0}^{k-1} \|\bx^{i} - \bx^{*} \| - \|\bx^{i+1} - \bx^{*} \| \nonumber \\
& \leq \|\bx^0 - \bx^{*} \| - \| \bx^k - \bx^0 \| \nonumber  \\
& \leq \| \bx^0 - \bx^{*}\|. \label{dppmtargetprovide3}
\end{align}
The derivation of the second inequality is based on the fact that the DPPM is a descent algorithm (Lemma \ref{dppmlem}). The derivation of the third inequality is due to (\ref{DPPMvalue1}).
Based on (\ref{dppmtargetprovide3}), the number of iterations $k$ required to attain an $\epsilon$-suboptimal solution for $f$ is of order $\frac{1}{\epsilon}$. This can be accelerated to ${\cal O}(\frac{1}{\sqrt \epsilon})$ using the Nesterov's acceleration method \cite{nesterov1983method}, as shown in Appendix~\ref{sec:acc}.

\qed

\subsection{Strongly convex quadratic functions}

Here, we characterize the rate of convergence for strongly convex quadratic functions. We expect that asymptotic  convergence rates obtained for quadratic functions are directly analogous to the general case, due to the fact that a twice continuously differentiable function $f: \R^n \rightarrow \R$ can be approximated near the local minimum $\bx^*$ using the following quadratic function:
\begin{align}\label{quadraticfunc}
f(\bx^*) + \frac{1}{2} (\bx - \bx^*)^\top \nabla^2 f(\bx^*) (\bx - \bx^*),
\end{align}
where $\nabla^2 f(\bx^*)$ is the Hessian matrix of $f$, which is a positive definite symmetric function.
Without a loss of the generality, $f$ can be minimized at $\bx^* = 0$ and $f(\bx^*) = 0$. Thus,
\begin{align*}
f(\bx) = \frac{1}{2}\bx^\top Q \bx, \text{     } \nabla f(\bx) = Q\bx, \text{ and    } \nabla^2 f(\bx) = Q.
\end{align*}
The DPP-update takes the following form:
\begin{align} \label{quadratic1}
\bx^{k+1} = \bx^k + w_k \bar \bp_k = \bx^k - t_k (\bar \bp_k^\top Q \bx^{k+1}) \bar \bp_k.
\end{align}
The second equality uses $w_k = - t_k \bar \bp_k^\top Q \bx^{k+1}$ from (\ref{dppmwvalue}). Rearranging (\ref{quadratic1}) yields
\begin{align} \label{quadupdate}
[I + t_k \bar \bp_k (\bar \bp_k^\top Q)] \bx^{k+1} = \bx^k.
\end{align}
Below, we show that $\bar \bp$ is the eigenvector of $[I + t\bar{\bp} \bar{\bp}^{\top} Q]$ with corresponding eigenvalue of $1+ t \bar \bp^\top Q \bar \bp$ for any $t > 0$.
\begin{citedlem}\label{PROPOSITION_1}
For any symmetric matrix $Q \in \mathbb{R}^{n\times n}$ and $t \in \mathbb R$, 
if $(1+t \bar{\bp}^{\top} Q \bar{\bp} )\neq 0$
then \\
(i) 
\begin{equation}\label{Prop_inv}
[\bI + t\bar{\bp} \bar{\bp}^{\top} Q]^{-1} 
= \bI - \frac{t}{(1+t \bar{\bp}^{\top} Q \bar{\bp} )}\bar{\bp} \bar{\bp}^{\top} Q,
\end{equation}
where $\bI$ is the $n\times n$ dimensional identity matrix. \\
(ii) $\bar \bp$ is the eigenvector of $[\bI + t\bar{\bp} \bar{\bp}^{\top} Q]^{-1}$ with corresponding eigenvalue of $\frac{1}{1+ t \bar \bp^\top Q \bar \bp}$. The eigenvalues of the other eigenvectors of $[\bI + t\bar{\bp} \bar{\bp}^{\top} Q]^{-1}$ are all equal to $1$.
\end{citedlem}
\proof

See Appendix \ref{inversequad}.

\qed

In accordance with (\ref{quadupdate}), $Q$ is a positive definite symmetric matrix, and based on Lemma \ref{PROPOSITION_1}(i), we obtain
\begin{align}\label{it_relation}
\bx^{k+1}  = (\bI - \frac{t}{(1+t \bar{\bp}^{\top} Q \bar{\bp} )}\bar{\bp} \bar{\bp}^{\top} Q) \bx^k = [\bI + t\bar{\bp} \bar{\bp}^{\top} Q]^{-1}  \bx^k .
\end{align}
Below, we demonstrate that it is possible to achieve an R-superlinear convergence rate using the DPPM. Let $\{\bar \bg_1, \cdots \bar \bg_n\}$ denote the conjugate vectors of unit 2-norm (that is, $\|\bg_i\| = 1$) with respect to $Q \in \R^{n \times n}$ if they are linear independent and $\bg_i^\top Q \bg_j = 0$ for $i \neq j$. Note that eigenvectors of $Q$ are also conjugate vectors. 
We can choose the cyclic descent-conjugate-vector as the search direction for the DPPM with all $l = 0, 1, \cdots$
\begin{align} \label{descenteigen}
\{ \bar  \bp_{ln+1} = \pm \bar \bg_{1}, \bp_{ln+2}= \pm \bar \bg_{2},  \cdots,  \bp_{(l+1)n} = \pm \bar \bg_{n}\}.
\end{align}
Then, the search direction of the cyclic descent-conjugate-vector is defined as follows:
\begin{align} \label{eigendir}
\bar \bp_{ln + i} = 
\begin{cases} 
\bar \bg_i \text{ if $\bg_i^\top Q \bx^{ln+i} < 0$,} \\
- \bar \bg_i \text{ otherwise. }
\end{cases}
\end{align}
Moreover, if the latter case in (\ref{eigendir}) is $\bg_i^\top Q \bx^{ln+i} =0$ (that is, neither $\bar \bg_i$ nor $-\bar \bg_i$ is a descent direction), then we let $\bx^{ln+(i+1)} = \bx^{ln+i}$.
In accordance with (\ref{descenteigen}), we can express $\bx^k$ as $\sum_{i = 0}^{n-1} b_i \bar{\bp}_{k+i}$.
Based on (\ref{it_relation}) and Lemma \ref{PROPOSITION_1}(ii), we have
\begin{align}{\label{dim_it_relation}}
    \bx^{k+n}
    &=(\bI+t_{k+n-1}\bar{\bp}_{k+n-1}  \bar{\bp}_{k+n-1}^{\top} Q )^{-1}
\cdots  
    (\bI+t_{k}\bar{\bp}_{k}  \bar{\bp}_{k}^{\top} Q )^{-1}\bx^k \nonumber \\
    &=
    \sum_{i=0}^{dim-1}\frac{b_i}{1 + t_{k+i} \bar{\bp}_{k+i}^\top Q \bar{\bp}_{k+i}}\bar{\bp}_{k+i} .
\end{align}
Suppose that $m$ is the smallest eigenvalue of $Q$ and $t_k = \lambda > 0$ for all $k$, due to the fact that the quadratic function is convex and therefore satisfies Lemma \ref{convexpara} immediately. From (\ref{dim_it_relation}), we obtain 
\begin{align} \label{quadupdatenorm*}
    \| \bx^{k+n} \|_Q     =     \| \sum_{i=0}^{n-1} \frac{b_i }{1 + \lambda \bar{\bp}_{k+i}^\top Q \bar{\bp}_{k+i}} \bar{\bp}_{k+i}\|_Q      \leq
    \|\sum_{i=0}^{n-1} \frac{b_i }{1 + \lambda m} \bar{\bp}_{k+i}\|_Q 
    =    \frac{1}{1 + \lambda m}  \|\bx^k\|_Q,
\end{align}
where $\|\cdot\|_Q$ denote the norm with respect to $Q$.

The convergence rate of the DPPM in the search direction of the cyclic descent-conjugate-vector is given as follows:
\begin{citedprop} \label{quadraticconvergence}
Suppose that $m$ is the smallest eigenvalue of $Q$. 
Let $\{ \bx^k \}$ be the sequence of iterates of the DPPM for $f$ (\ref{quadraticfunc}), the derivation of which is based on the search direction of the cyclic descent-conjugate-vector, defined in (\ref{eigendir}). Then, \\
(i) If we let $t_k = \lambda  > 0$ for all $k$, then sequence $\{ \bx^k \}$ converges to the zero vector (the optimal point $x*$) and the convergence rate is R-linear and equal to $(\frac{1}{1+\lambda m})^\frac{1}{n}$. \\
(ii) If $\lambda$ is a positive monotonic increasing function for which the number of of iterations satisfies     \begin{align}\label{increasing_lambda}
        \frac{\lambda((k+1)n)}{\lambda(kn)}
        &\geq
        c > 1,
    \end{align}
 then sequence $\{ \bx^k \}$ has 
    \begin{align*}
        \frac{\|\bx^{(k+1)n}\|_Q}{\|\bx^{kn}\|_Q}
        &\leq 
        \frac{1}{1 + c^k \lambda(0) m}, \text{ $\forall k \in \mathbb{N}\cup \{0\}$}.
    \end{align*}
    Hence, sequence $\{ \bx^k \}$ converges to optimal point $x^* = 0$, and the convergence rate is R-superlinear.
\end{citedprop}

\proof
(i) From (\ref{quadupdatenorm*}), we can obtain for any $k$
\begin{align*}
\frac{\| \bx^{kn} \|_Q}{\|\bx^0\|_Q} \leq [\frac{1}{1 + \lambda m}]^k.
\end{align*}
From $0<\frac{1}{1+\lambda m}<1$, we obtain
\begin{align*}
    \lim_{k \rightarrow \infty} \frac{\|\bx^{kn}\|_Q}{\|\bx^0\|_Q}
    &=
    \lim_{k \rightarrow \infty} [\frac{1}{1+\lambda m}]^{k} = 0.
\end{align*}
Thus, $\{\|\bx^{kn}\|_Q\}_k$ converges to zero as $k$ tends toward infinity. We also know that sequence 
$\{\|\bx^k\|_Q\}$ of the DPPM is a monotone non-increasing sequence. Thus, we can deduce that $\{\|\bx^k\|_Q\}$ converges to zero, which implies that sequence $\{ \bx^k\}$ converges to the zero vector.

Let the error sequence be $\{\epsilon_k\}$, where $\epsilon_k = \|\bx^k - \bx^*\|_Q= \|\bx^k\|_Q$. In accordance with (\ref{quadupdatenorm*}), $\{\epsilon_k\}$ satisfies
\begin{align*}
        \frac{\epsilon_{(k+1)n}}{\epsilon_{kn}}
        &\leq 
        \frac{1}{1 + \lambda m}, \text{$\forall k \in \mathbb{N}\cup \{0\}$}.
\end{align*}
As $\{\epsilon_k\}$ is a monotone non-increasing sequence, in accordance with the definition of R-linear, we can conclude that the convergence rate of the DPPM is R-linear and equal to $(\frac{1}{1+\lambda m})^\frac{1}{n}$ for $\lambda > 0$.

(ii) Based on the fact that $\lambda$ is a positive monotonic increasing function for which the number of iterations satisfies (\ref{increasing_lambda}), we obtain
\begin{align*}
        \frac{\epsilon_{(k+1)n}}{\epsilon_{kn}}
        \leq 
        \frac{1}{1 + \lambda(kn) \cdot m}
        \leq 
        \frac{1}{1 + c^k \lambda(0)  m}, \text{ $\forall k \in \mathbb{N}\cup \{0\}$}.
\end{align*}
Because $c>1$, we have
\begin{align*}
        \limsup_{k \rightarrow \infty}\frac{\epsilon_{(k+1)n}}{\epsilon_{kn}}
        \leq 
        \limsup_{k \rightarrow \infty}\frac{1}{1 + c^k \lambda(0)  m}
        =
        0.
\end{align*}
Using an argument similar to the proof in (i), this inequality demonstrates that sequence $\{ \bx^k \}$ converges to optimal point $x^* = 0$, and that the convergence rate is R-superlinear.

\qed


\section{Algorithm and experiments}\label{experiment}

Algorithm \ref{DPPMalg} outlines the DPPM in our experiments. Steps 1 and 3 determine the descent direction and the convex segment for the DPPM in each iteration.

\begin{algorithm}[!h] 
\caption{The DPPM } 
 \label{DPPMalg}
\algsetup{indent=2em}
\begin{algorithmic}[1]
\item[\bf Input:]  function $f$ and the initial iterate $\bx^0$. 
\STATE Select a search direction.
\STATE Determine the convex segment and the optimal step-size. 
\STATE Update to the next iterate.
\STATE If a stopping condition is not reached, then go to Step 1.
\end{algorithmic}
\end{algorithm}

Algorithm \ref{dppmalgo} describes Step 2 of Algorithm \ref{DPPMalg}, which  determines the convex segment $[0, v(\bx)]$ at the current iterate $\bx$ (when direction $\bar \bp$ is provided) and provides the solution to (\ref{dppmwform}) with optimal step-size $w^*$, thereby satisfying 
\begin{align*} 
0= w + t \bar \bp^\top \nabla f(\bx + w \bar \bp),
\end{align*}
where parameter $t$ is derived from the convex segment.
The algorithm assumes that $(0, v(\bx, \bar \bp)] \subseteq (0, 1]$ is the convex lower-level set of $f$ at $\bx$ along ray $\bar \bp$, using $v(\bx, \bar \bp)$ derived in Steps 1-3. 
In Step 4, we use the golden section search method \cite{avriel1968golden} to find the unique minimum of 
(\ref{dppmwform})(Lemma \ref{dppmunique}(ii)). The fact that the golden search method does not use any  gradient evaluations whatsoever makes it ideally suited to situations in which the gradient of a function cannot be efficiently or accurately derived.
A minimum is known to be bracketed when there is a triplet point, $a < b < c$, such that $f(\bx + b \bar \bp)$ is larger than $f(\bx + a \bar \bp)$ and $f(\bx + c \bar \bp)$. This method involves selecting a new point $z$, either between $a$ and $b$ or between $b$ and $c$. If we make the latter choice, then we evaluate $f(\bx + z \bar \bp)$. If $f(\bx + b \bar \bp) > f(\bx + z \bar \bp)$, then the new bracketing triplet of points is $(a, b, z)$; otherwise, if $f(\bx + b \bar \bp) < f(\bx + z \bar \bp)$, then the new bracketing triplet is $(b, z, c)$. Parameter $t$ in Step 3 is derived from $v(\bx, \bar \bp)$ in accordance with Lemma \ref{dppmunique} (iii) in Step 4 to guarantee that $w^*$ falls within the interval. Furthermore, if $f$ is a convex function over $\Omega_l=\{\bx | f(\bx) \leq f(\bx^l)\}$ for a given $l$, then in accordance with Corollary \ref{convexpara}, $t_k$ can be a constant for $k > l$. In this case, Step 2 of Algorithm  \ref{dppmalgo} can be skipped.

\begin{algorithm}[t] 
\caption{Step 2 of Algorithm \ref{DPPMalg}} 
 \label{dppmalgo}
\algsetup{indent=2em}
\begin{algorithmic}[1]
\item[\bf Input:]  $\bar \bp$ is a direction at $\bx$; $\tau > 0$ is a parameter. 
\item[\bf Output:] $w^* \in (0, v(\bx, \bar \bp)]$, where $(0, v(\bx, \bar \bp)] \subseteq (0, 1]$ is the convex lower level interval of $f$ at $\bx$ along the ray $\bar \bp$.
\STATE Partition the interval $[\bx, \bx+ \bar \bp]$ into $10^{\tau}$ pieces of equal length segments where end points are at $\bx$, $\bx+\frac{1}{10^\tau} \bar \bp$, $\bx+\frac{2}{10^\tau} \bar \bp$, $\cdots$, $\bx + \bar \bp$ and evaluate the function values at the points to obtain $f(\bx)$, $f(\bx+\frac{1}{10^\tau} \bar \bp)$, $f(\bx+\frac{2}{10^\tau} \bar \bp)$, $\cdots$, $f(\bx + \bar \bp)$.
\STATE Determine the convex lower level set $[\bx, \bx + v(\bx, \bar \bp) \bar \bp] \subseteq [\bx, \bx+ \bar \bp]$. \\
If there exists $10^\tau-1\ge k \ge 0$, $k$ is the smallest integer such that $f(\bx+\frac{k}{10^\tau} \bar \bp) - f(\bx+\frac{k-1}{10^\tau} \bar \bp) > f(\bx+\frac{k+1}{10^\tau} \bar \bp) - f(\bx+\frac{k}{10^\tau} \bar \bp)$, then $v(\bx, \bar \bp) = \bx+\frac{k}{10^\tau} \bar \bp$; otherwise, $v(\bx, \bar \bp) = 1$.
\STATE Set the value of parameter $t$ in accordance with Lemma \ref{dppmunique}.
\STATE Apply the golden section search method \cite{press1988numerical} to obtain $w^*$.
\RETURN $w^*$.
\end{algorithmic}
\end{algorithm}

We conducted two experiments involving the DPPM. The first experiment demonstrated that the DPPM can converge to a critical point of a non-convex function. The second experiment was meant to confirm that the convergence rate for strongly convex quadratic functions is indeed R-superlinear.

\begin{figure}[H]
    \centering
    \includegraphics[scale = 0.7]{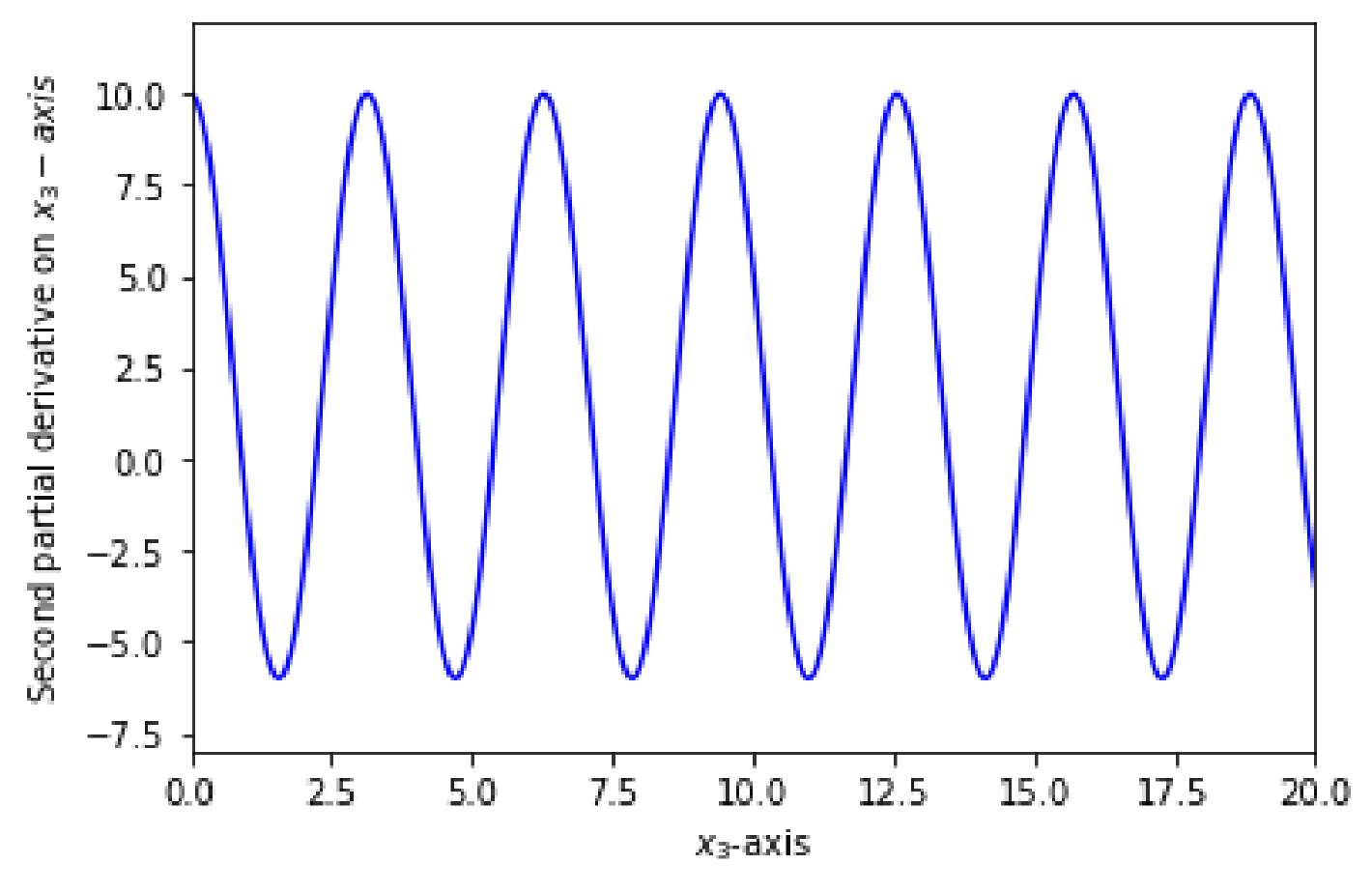}
    \caption{Plot of second partial derivative of $f(\bx = (x_1, x_2, x_3)) = \|\bx\|^2 + 4\cdot \sin^2{x_3}$ on $x_3$-axis where regions with $\frac{\partial^2}{{\partial x_3}^2} f < 0$ are non-convex.}
    \label{fig:nonconvex region}
\end{figure}

\begin{figure}[H]
\centering
\subfigure[]
{  \includegraphics[scale = 0.5]{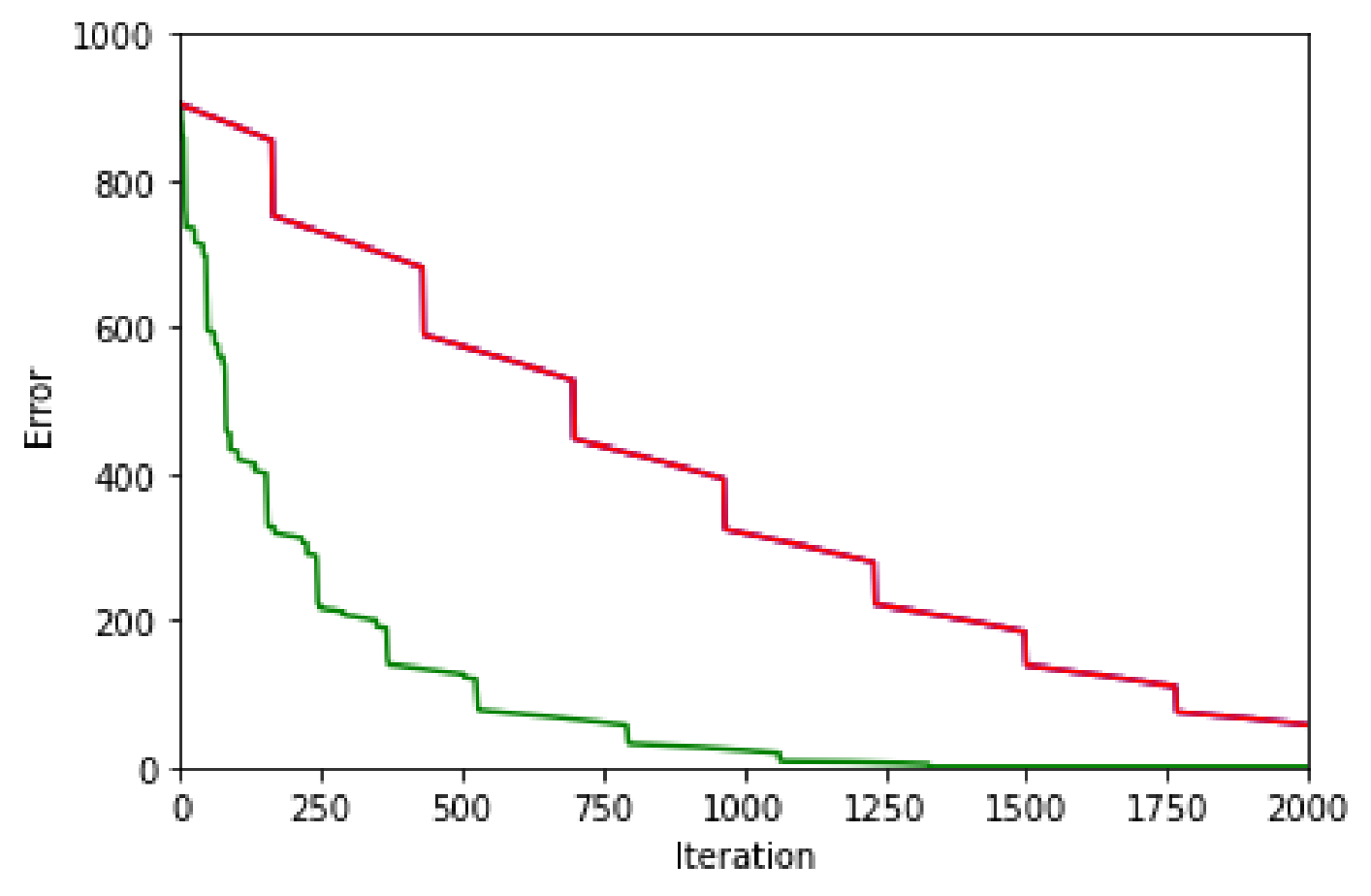}}
\subfigure[]
{  \includegraphics[scale = 0.5]{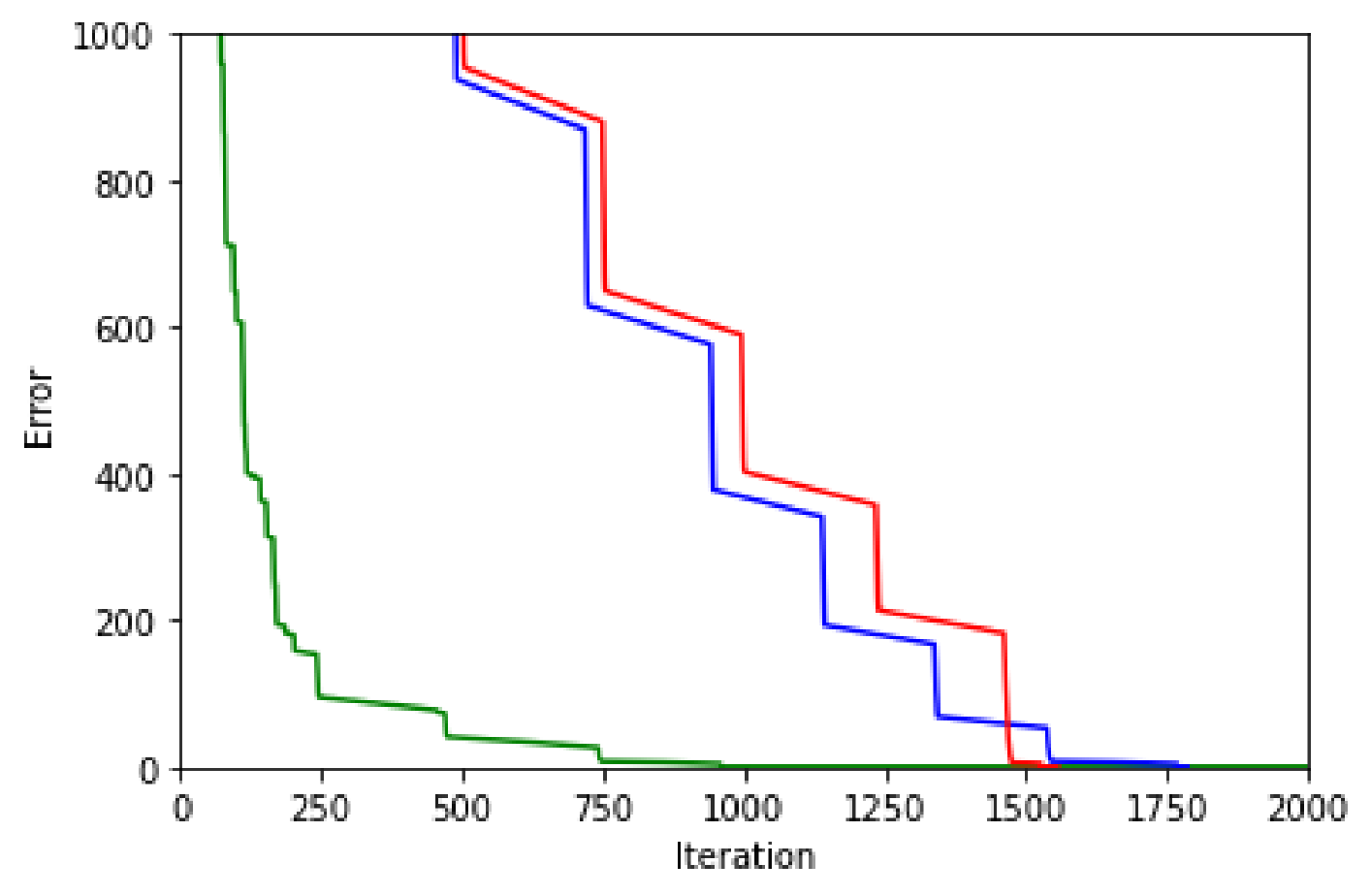}}
\subfigure[]
{  \includegraphics[scale = 0.5]{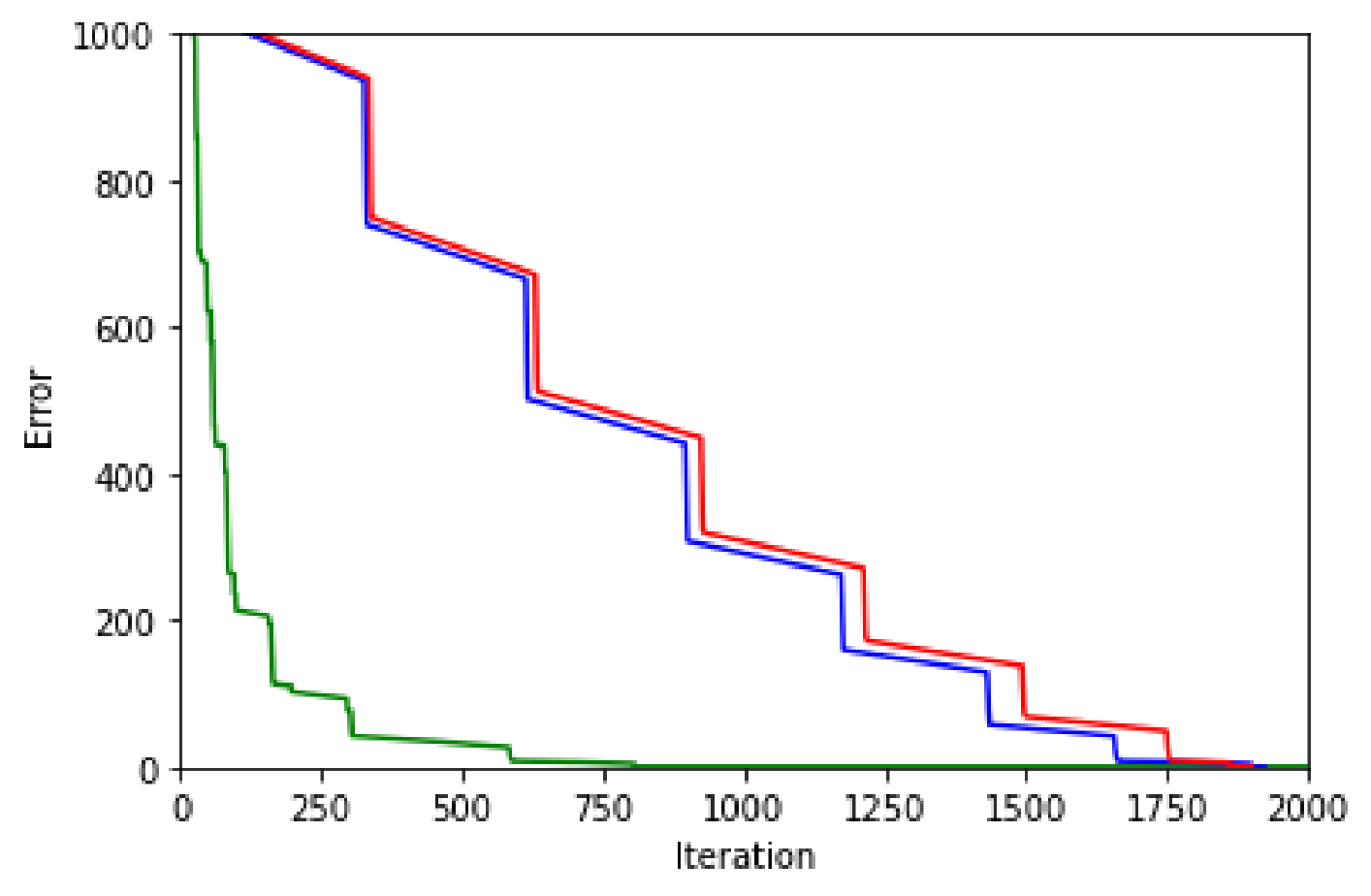}}
\subfigure[]
{  \includegraphics[scale = 0.5]{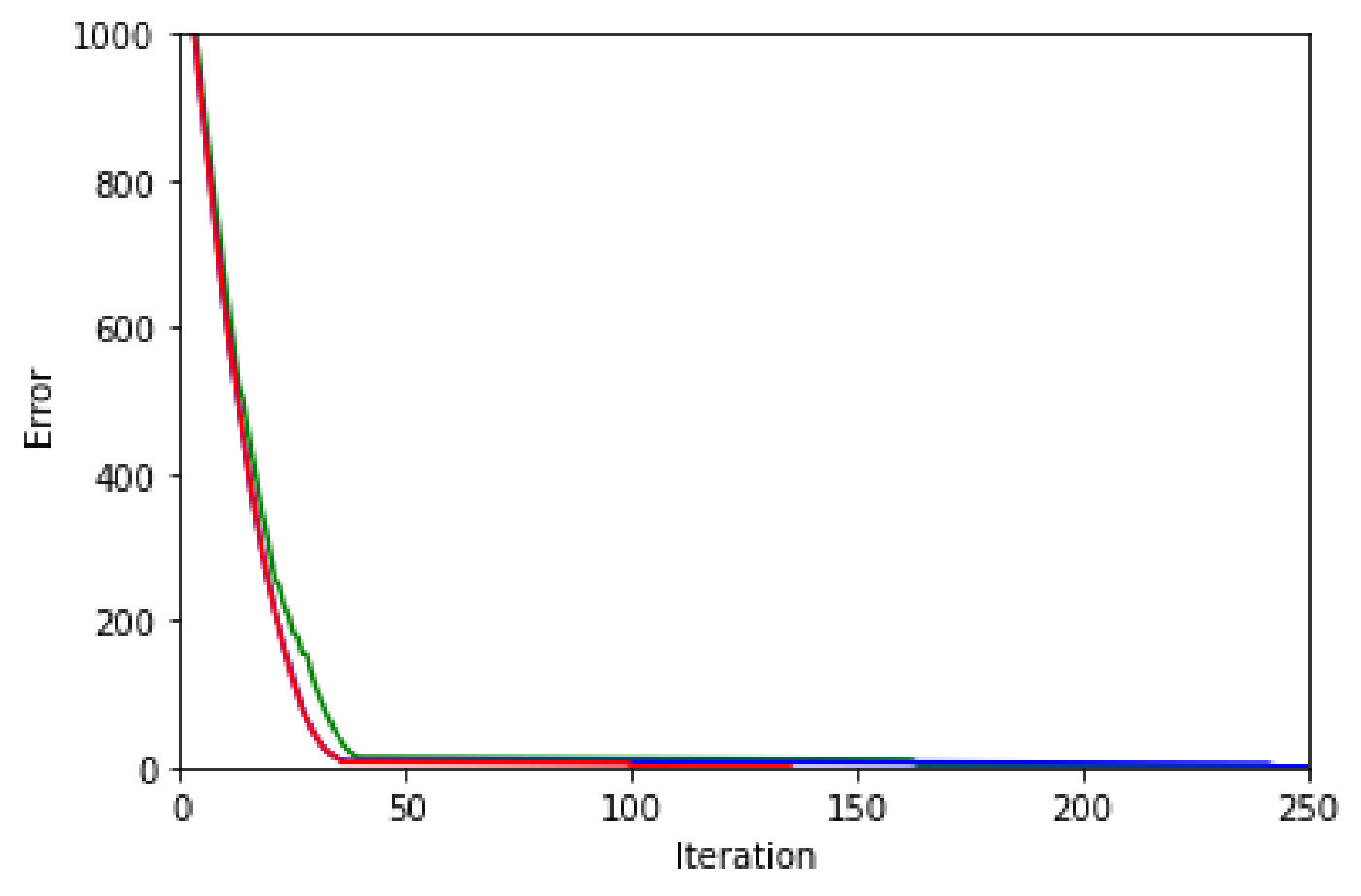}}
    \caption{Plot of error vs. iterations where blue, red, and green lines respectively indicate results obtained using gradient, momentum, and DLC direction (Section \ref{sec:DLCdirection}) methods to determine search direction $\bar{\bp}_k$.  (a) The initial guess selected at the point $[0, 0, 30]$ on the $x_3$-axis. The blue and red lines are overlapped. (b)-(d) The initial guess selected as random with $10 \leq x_i^{(0)} \leq 40$ for $i= 1, 2, 3$. Momentum weight set at $\beta_k=0.6$ for all iterations \cite{liu2020improved}. Parameter $\mu = 1000$ and $\bp_0$ was selected from $\nabla f(\bx) + j*0.5* |\nabla f(\bx)|* q$ where $q$ is a random unit norm vector and $j$ is from $0$ to $9$ at random for the DLC direction method. In all cases, the DLC direction method is superior to the compared methods.}
    \label{fig:nonconvex}
\end{figure}

In the first experiment, we considered the non-convex function $f(\bx = (x_1, x_2, x_3)) = \|\bx\|^2 + 4\cdot \sin^2{x_3}$. As shown in Figure \ref{fig:nonconvex region}, this function shows non-convexity along the $x_3$ direction, due to the negative second derivative on some of the segments along the axis. 
Figure \ref{fig:nonconvex} compares the convergence of  $\{f(\bx^k)\}$ to the optimal value $f^*$ of the function using the directions of descent, derived using the gradient (blue), momentum (red), and DLC direction (green) methods.

Note that for this function, if an iterate lies on the $x_3$-axis (corresponding to the negative second derivative regions in Figure \ref{fig:nonconvex region}), then the gradient direction does not satisfy (DCL). In this case, the gradient direction is not legitimate for DPPM. Thus, we imposed a slight perturbation in this iterate to make the gradient direction legitimate for the DPPM.

\begin{figure}[H]
\centering
\subfigure[$\lambda = 0.1$.]{
\begin{minipage}[t]{0.5\linewidth}
\centering
\includegraphics[width=7cm]{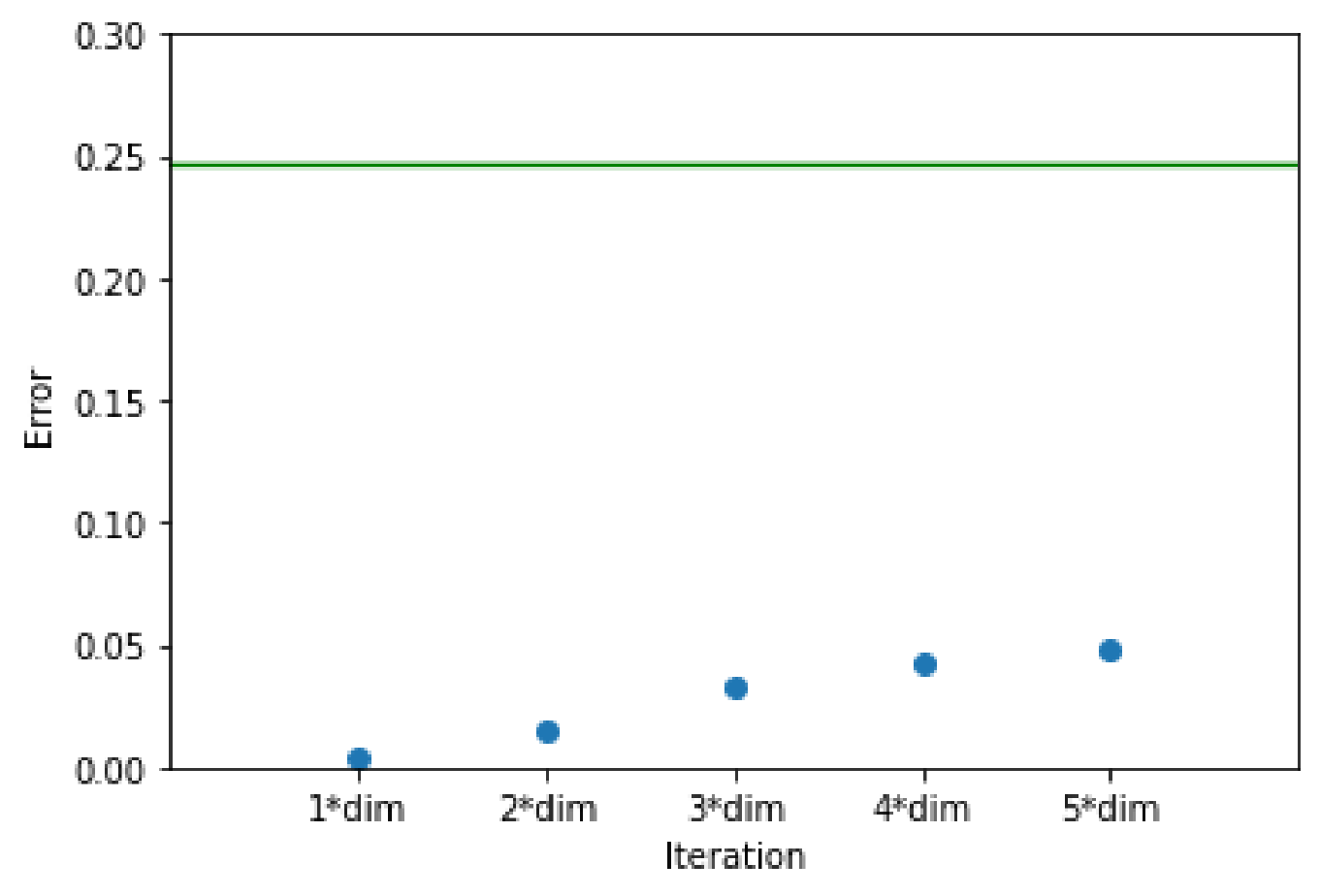}
\end{minipage}%
}%
\subfigure[$\lambda = 10$.]{
\begin{minipage}[t]{0.5\linewidth}
\centering
\includegraphics[width=7cm]{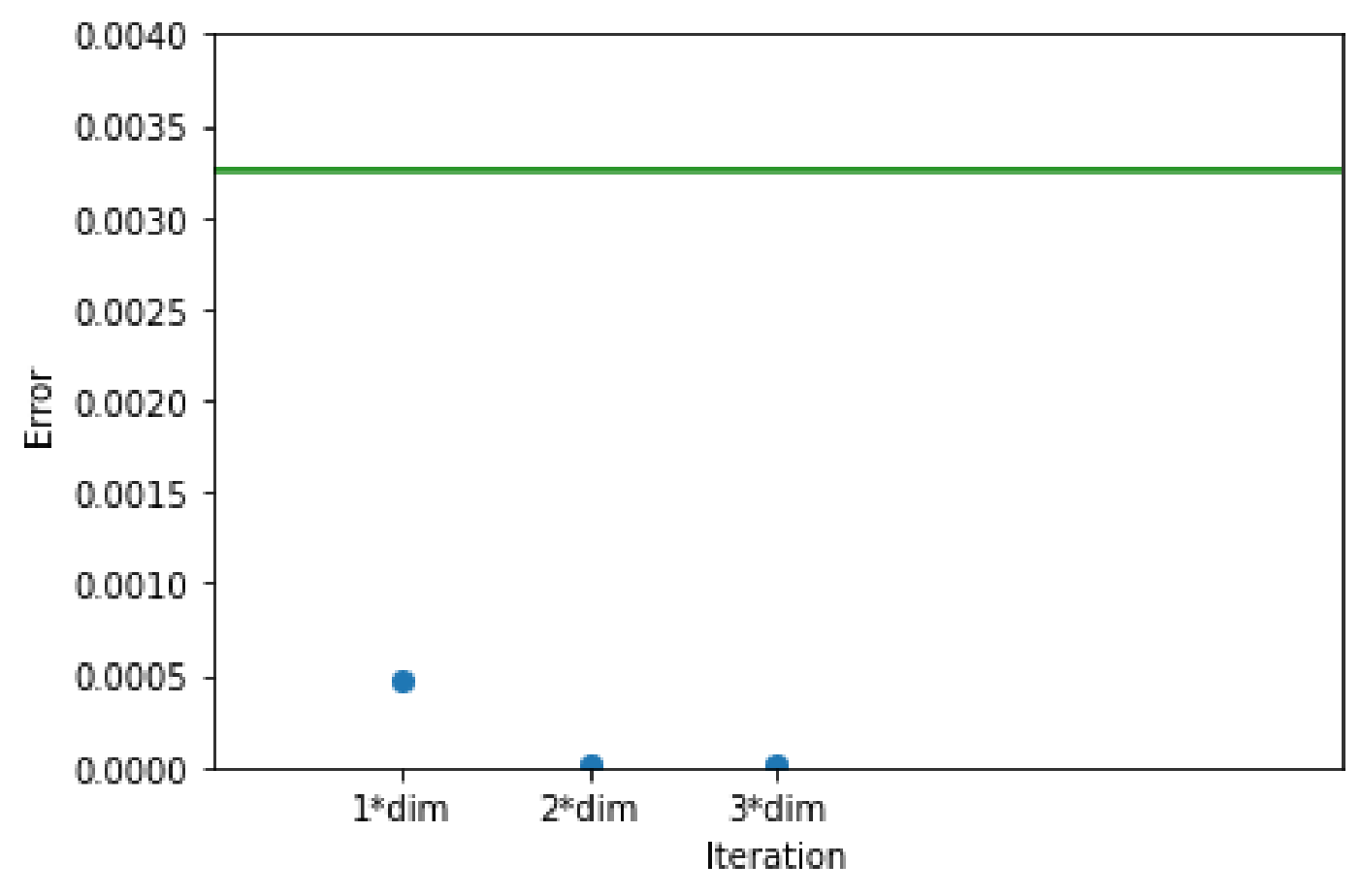}
\end{minipage}%
}%
\centering
\caption{Plot of ratio vs. iterations for different values of $\lambda$. In each sub-figure, blue dots denote ratio $\frac{\|\bx^{(k+1)\cdot dim}\|_Q}{\|\bx^{k\cdot dim}\|_Q}$ with  $dim = 500$ and green line indicates theoretical upper bound ($\frac{1}{1+\lambda m}$, where $m$ is the smallest eigenvalue of the strongly convex quadratic function). }\label{fig:eigenvector}
\end{figure}

The second experiment deals with the convergence rate of a strongly convex quadratic function, which is $\frac{1}{2} \bx^\top Q \bx$, where $Q \in \mathbb R^{500 \times 500}$ is a diagonal matrix with diagonal elements generated from a $[0, 1]$ uniform random variable. The value of a diagonal  element is then scaled to $[30, 300]$. The optimal solution of this function is at the origin. In accordance with Corollary \ref{convexpara}, we set parameter $t = \lambda$ as a constant and used the cyclic descent-conjugate-vector search direction (\ref{eigendir}) in the search for the descent direction. The conjugate vectors adopted in our experiment were eigenvectors of $Q$. Because $Q$ is a diagonal matrix, the eigenvectors of $Q$ are standard basis.
As shown in Figure \ref{fig:eigenvector}, the ratio of the error sequence fell under the theoretical bound given in Proposition \ref{quadraticconvergence}(i) for $\lambda = 0.1$ (a) and $10$ (b). 
This indicates that the DPPM achieved an $R$-linear convergence rate for the strongly convex quadratic function with constant $\lambda$ values. 

\begin{figure}[tp]
\centering
\subfigure[Green line is $\frac{1}{1+\lambda(0) \cdot m}$.]{
\includegraphics[width=7.7cm]{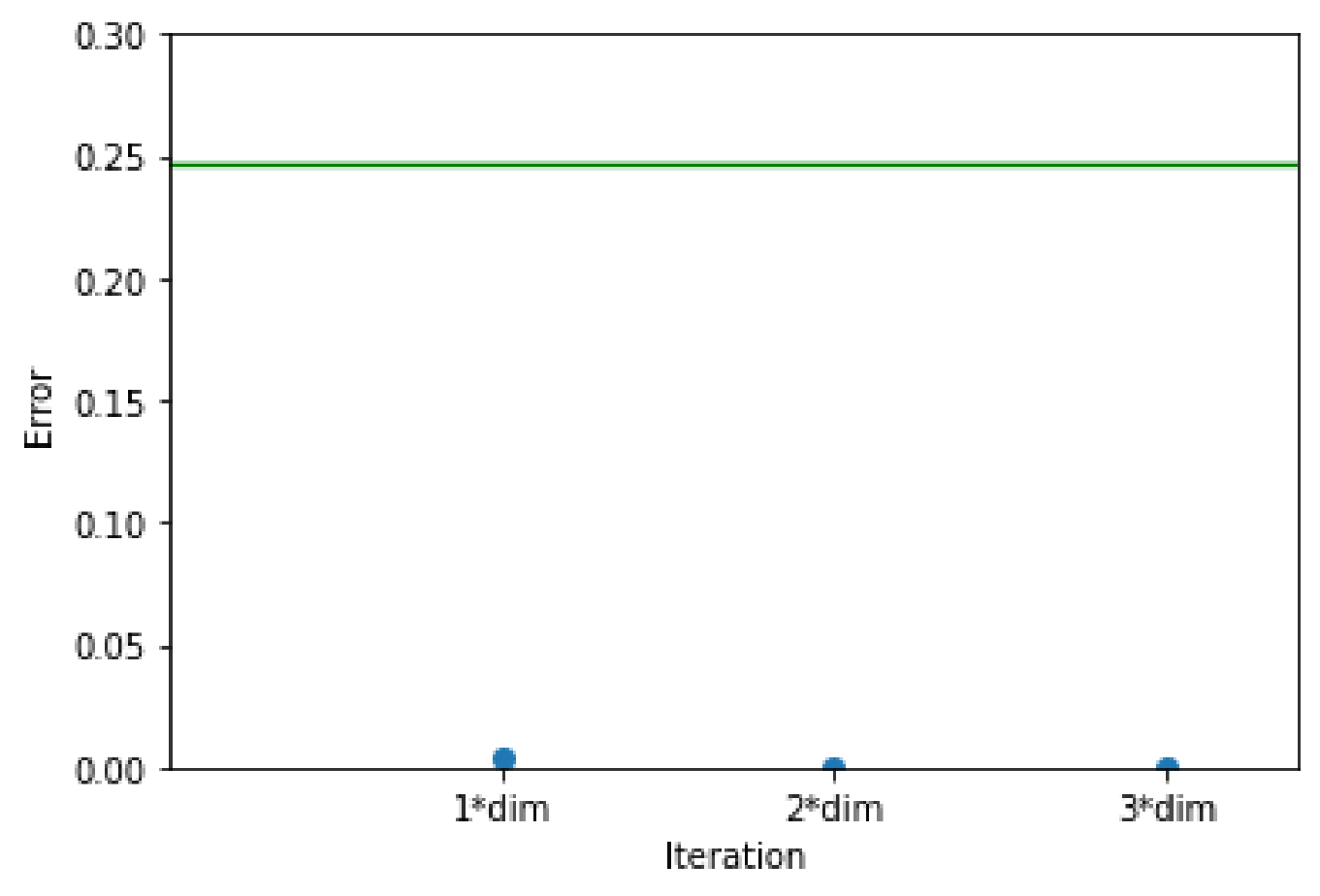}
}
\subfigure[Yellow line is $\frac{1}{1+c\lambda(0) \cdot m}$.]{
\includegraphics[width=7.7cm]{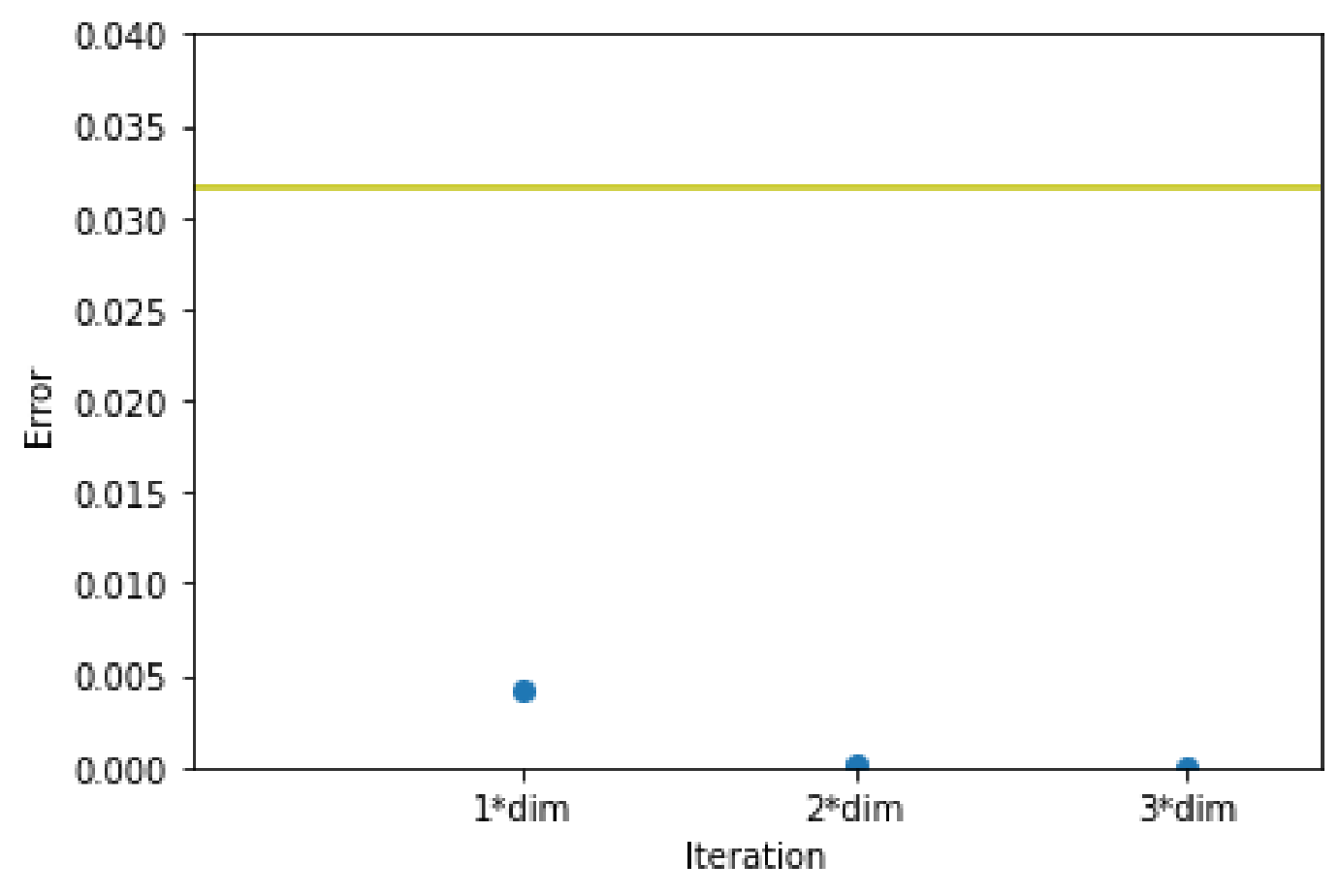}
}
\subfigure[Red line is $\frac{1}{1+c^2\lambda(0) \cdot m}$.]{
\includegraphics[width=7.7cm]{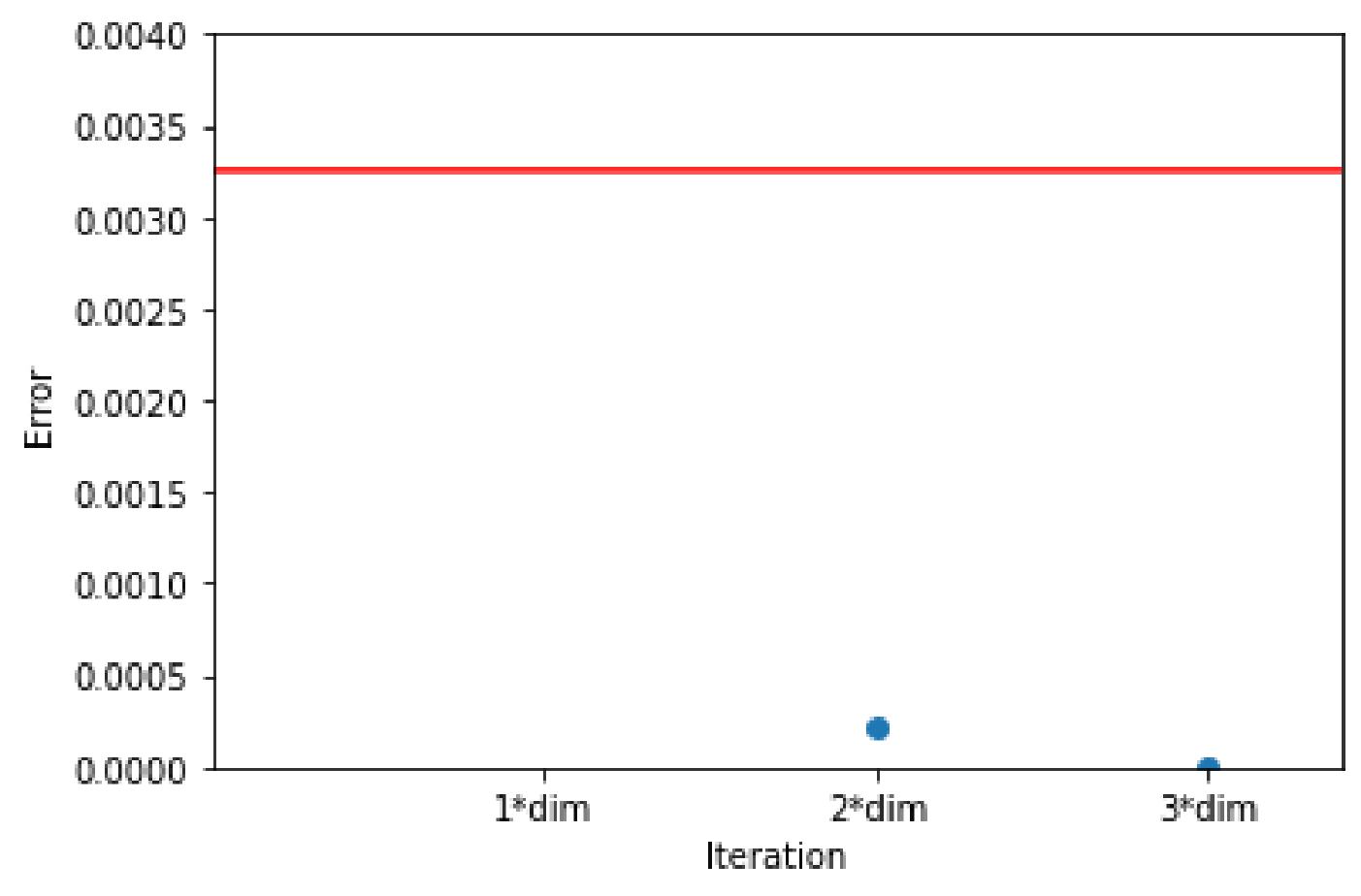}
}

\caption{Plot of ratio vs. iterations where blue dots in (a), (b), and (c) denote the ratios $\frac{\|\bx^{(k+1)\cdot dim}\|}{\|\bx^{k\cdot dim}\|}$ where $dim=500$ is the number of variables for $k=1$ ($1*$dim), $k=2$ ($2*$dim), and $k=3$ ($3*$dim), and green, yellow, and red lines are respectively $\frac{1}{1+\lambda(0) \cdot m}$, $\frac{1}{1+c\cdot\lambda(0) \cdot m}$, and $\frac{1}{1+c^2\cdot\lambda(0) \cdot m}$, where $m$ is the minimum eigenvalue of the quadratic function. Parameter $t$ was set to be monotonically non-decreasing as the number of iterations increases. This ensures the error sequence converges super-linearly to zero. For the first $500$ DPPM iterates, parameter $t$ was set at $\lambda(0)$;  for $501$ to $1000$, it was set at $\lambda(0) c$; and it was set at $\lambda(0) c^2$ for the last $500$ iterates, where $c =10$ and $\lambda(0) = 0.1$. Due to the wide range of the ratio values, (a)-(c) were plotted to compare the decreasing error ratios to the analytical bounds, where each bound is obtained with a constant parameter $t$. In (a), all three blue dots ($k=1,2, 3$) fall below the bound; in (b), the last two blue dots ($k=2, 3$) fall below the bound; and in (c), the last blue dots ($k=3$) fall below the line. No bound fell below the blue dots - an indication of super-linear convergence to zero.}\label{eigenvector_increasing}
\end{figure}

To achieve the R-super-linearity, implied by Lemma \ref{quadraticconvergence} (iv), we increased the value of $\lambda$ (and accordingly $t$) with the number of iterations. As shown in Figure \ref{eigenvector_increasing}, we then demonstrated that  $\{\bx^k\}$ converges to the optimal point as $\frac{\|\bx^{(k+1)\cdot dim}\|_Q}{\|\bx^{k\cdot dim}\|_Q} \leq \frac{1}{1 + c^k \cdot\lambda(0)\cdot m} \rightarrow 0$ as $k \rightarrow \infty$, where $dim$ indicates the dimension of the variables, $c=10$, and $\lambda(0) =0.1$.

\section{Conclusions} \label{sec:conclusion}

This paper introduces the directional proximal point method (DPPM) by which to solve the problem of un-constrained minimization for smooth but not strictly concave functions. We demonstrate that the search direction and optimal step-size for an update of an iterate can both be derived by optimizing no more than two variables, regardless of the dimensionality of the function variables.  
This gives the DPPM an advantage over comparable methods when dealing with large-scale problems.
We demonstrate that if the sequence of the directions is gradient-related, then the DPPM can converge to critical points. We also present conditions pertaining to the descent directions and local properties of a critical point for the entire DPPM sequence to converge to a single critical point. When dealing with convex functions, the DPPM is as efficient as the PPM and can be accelerated using the Nestorov approach. For strongly convex quadratic functions, the error sequence of DPPM converges R-superlinerly to zero, regardless of the dimensionality of the function variables.
The DPPM could potentially be extended to include the problem in which $f$ satisfies (DLC) but $f$ is not smooth.  
This extension would move a step closer to solving the more general problem of
\begin{align*}
\min_{\bx \in \chi} f(\bx) = g(\bx) + h(\bx),
\end{align*}
where $g$ and $h$ respectively refer to smooth and non-smooth functions involving a large number of variables.

\section*{Acknowledgments}
Wen-Liang Hwang would like to thank Professor Andreas Heinecke at Yale-NUS college for his valuable comments. He would also like to express gratitude to the authors of  the lecture notes accessible to the public from which we have benefited greatly in terms of teaching and research. The public available software for the DPPM is given in https://github.com/Mick048/DPPM

%
%
%
%
%
%

\bibliographystyle{ieeetr}
\bibliography{DeepRef}

\begin{thebibliography}{10}

\bibitem{boyd2004convex}
S.~P. Boyd and L.~Vandenberghe, {\em Convex optimization}.
\newblock Cambridge university press, 2004.

\bibitem{bertsekas2014constrained}
D.~P. Bertsekas, {\em Constrained optimization and Lagrange multiplier
  methods}.
\newblock Academic press, 2014.

\bibitem{van2020}
L.~Vandenberghe, {\em ECE236C-Optimizatoin methods for large-scale systems}.
\newblock Lecture note of UCLA, 2020.

\bibitem{pil2020}
M.~Pilanci, {\em EE364b-Convex Optimizatoin II}.
\newblock Lecture note of Stanford, 2020.

\bibitem{ruder2016overview}
S.~Ruder, ``An overview of gradient descent optimization algorithms,'' {\em
  arXiv preprint arXiv:1609.04747}, 2016.

\bibitem{trefethen1997numerical}
L.~N. Trefethen and D.~Bau~III, {\em Numerical linear algebra}, vol.~50.
\newblock Siam, 1997.

\bibitem{bertsekas2015convex}
D.~P. Bertsekas and A.~Scientific, {\em Convex optimization algorithms}.
\newblock Athena Scientific Belmont, 2015.

\bibitem{bauschke2011convex}
H.~H. Bauschke, P.~L. Combettes, {\em et~al.}, {\em Convex analysis and
  monotone operator theory in Hilbert spaces}, vol.~408.
\newblock Springer, 2011.

\bibitem{Ber08}
D.~P. Bertsekas, ``Nonlinear programming: 3rd,'' {\em Athena Scientific
  Optimization and Computations Series 4}, vol.~4, 2008.

\bibitem{nocedal2006numerical}
J.~Nocedal and S.~Wright, {\em Numerical optimization}.
\newblock Springer Science \& Business Media, 2006.

\bibitem{barzilai1988two}
J.~Barzilai and J.~M. Borwein, ``Two-point step size gradient methods,'' {\em
  IMA journal of numerical analysis}, vol.~8, no.~1, pp.~141--148, 1988.

\bibitem{raydan1993barzilai}
M.~Raydan, ``On the barzilai and borwein choice of steplength for the gradient
  method,'' {\em IMA Journal of Numerical Analysis}, vol.~13, no.~3,
  pp.~321--326, 1993.

\bibitem{dai2002r}
Y.-H. Dai and L.-Z. Liao, ``R-linear convergence of the barzilai and borwein
  gradient method,'' {\em IMA Journal of Numerical Analysis}, vol.~22, no.~1,
  pp.~1--10, 2002.

\bibitem{burdakov2019stabilized}
O.~Burdakov, Y.-H. Dai, and N.~Huang, ``Stabilized barzilai-borwein method,''
  {\em Journal of Computational Mathematics}, pp.~916--936, 2019.

\bibitem{fletcher2005barzilai}
R.~Fletcher, ``On the barzilai-borwein method,'' in {\em Optimization and
  control with applications}, pp.~235--256, Springer, 2005.

\bibitem{armijo1966minimization}
L.~Armijo, ``Minimization of functions having lipschitz continuous first
  partial derivatives,'' {\em Pacific Journal of mathematics}, vol.~16, no.~1,
  pp.~1--3, 1966.

\bibitem{fridovich2020approximately}
S.~Fridovich-Keil and B.~Recht, ``Approximately exact line search,'' {\em arXiv
  preprint arXiv:2011.04721}, 2020.

\bibitem{asl2020analysis}
A.~Asl and M.~L. Overton, ``Analysis of the gradient method with an
  armijo--wolfe line search on a class of non-smooth convex functions,'' {\em
  Optimization methods and software}, vol.~35, no.~2, pp.~223--242, 2020.

\bibitem{raydan1997barzilai}
M.~Raydan, ``The barzilai and borwein gradient method for the large scale
  unconstrained minimization problem,'' {\em SIAM Journal on Optimization},
  vol.~7, no.~1, pp.~26--33, 1997.

\bibitem{grippo1986nonmonotone}
L.~Grippo, F.~Lampariello, and S.~Lucidi, ``A nonmonotone line search technique
  for newton’s method,'' {\em SIAM journal on Numerical Analysis}, vol.~23,
  no.~4, pp.~707--716, 1986.

\bibitem{zhou2006gradient}
B.~Zhou, L.~Gao, and Y.-H. Dai, ``Gradient methods with adaptive step-sizes,''
  {\em Computational Optimization and Applications}, vol.~35, no.~1,
  pp.~69--86, 2006.

\bibitem{kinderlehrer2000introduction}
D.~Kinderlehrer and G.~Stampacchia, {\em An introduction to variational
  inequalities and their applications}.
\newblock SIAM, 2000.

\bibitem{pshenichnyj2012linearization}
B.~N. Pshenichnyj, {\em The linearization method for constrained optimization},
  vol.~22.
\newblock Springer Science \& Business Media, 2012.

\bibitem{bolte2014proximal}
J.~Bolte, S.~Sabach, and M.~Teboulle, ``Proximal alternating linearized
  minimization for nonconvex and nonsmooth problems,'' {\em Mathematical
  Programming}, vol.~146, no.~1, pp.~459--494, 2014.

\bibitem{nesterov1983method}
Y.~Nesterov, ``A method for unconstrained convex minimization problem with the
  rate of convergence o (1/k\^{} 2),'' in {\em Doklady an ussr}, vol.~269,
  pp.~543--547, 1983.

\bibitem{avriel1968golden}
M.~Avriel and D.~J. Wilde, ``Golden block search for the maximum of unimodal
  functions,'' {\em Management Science}, vol.~14, no.~5, pp.~307--319, 1968.

\bibitem{press1988numerical}
W.~H. Press, S.~A. Teukolsky, W.~T. Vetterling, and B.~P. Flannery, ``Numerical
  recipes in c,'' 1988.

\bibitem{liu2020improved}
Y.~Liu, Y.~Gao, and W.~Yin, ``An improved analysis of stochastic gradient
  descent with momentum,'' {\em NIPS}, 2020.

\bibitem{nesterov2013introductory}
Y.~Nesterov, {\em Introductory lectures on convex optimization: A basic
  course}, vol.~87.
\newblock Springer Science \& Business Media, 2013.

\bibitem{zangwill1969nonlinear}
W.~I. Zangwill, {\em Nonlinear programming: a unified approach}, vol.~52.
\newblock Prentice-hall Englewood Cliffs, NJ, 1969.

\end{thebibliography}

\appendix

\section{Smooth functions satisfying (DLC)} \label{AS}


\begin{citedlem}(Lemma 3.1.3 \cite{nesterov2013introductory})\label{strictlyconcave}
Let function $f$ be strictly concave. Then, for all $\bx, \bar \bp$ and $w \geq 0$, we have
\begin{align*}
f(\bx + w \bar \bp) < f(\bx) + \beta f'(\bx, \bar \bp)
\end{align*}
\end{citedlem}

\begin{citedlem} \label{strictlyconcave1}
Let $f \in C^1$ and for all $\bx$ and $\bp$,
\begin{align*}
f(\bx +\bp) < f(\bx) + \nabla f(\bx)^\top  \bp.
\end{align*}
Then $f$ is a strictly concave.
\end{citedlem}
\proof
%
%
%

\proof
Consider two distinct points $\bx_1$ and $\bx_2$ and let $\bx= \beta \bx_1 + (1-\beta) \bx_2$ where $\beta \in [0, 1]$. We can denote the $\bx_1$ and $\bx_2$ as:
\[\bx_1 = \bx + \bp_1,\]
\[\bx_2 = \bx + \bp_2.\]
Clearly, this also means that $\beta \bp_1 + (1-\beta) \bp_2 = 0$. Then, we have
\begin{align*}
    \beta f(\bx + \bp_1) + (1-\beta)  f(\bx + \bp_2)
    &<
    \beta (\nabla f(\bx)^\top \bp_1+ f(\bx)) + (1-\beta)  (\nabla f(\bx)^\top \bp_2 + f(\bx))\\
    &=
    f(\bx) + \nabla f(\bx)^\top (\beta \bp_1 + (1-\beta) \bp_2)\\
    &=
    f(\bx).
\end{align*}
This proof is done.

\qed


The class of function $f \in C^1$ satisfies (DLC) is characterized in the following corollary.

\begin{citedcor} \label{characterizeDLC}
Let $f \in C^1$. $f$ satisfies (DLC) if and only if $f$ is not a strict concave function.
\end{citedcor}
\proof

Assuming $f$ satisfies (DLC) at $\bx, \bar \bp$ and $w \geq 0$: 
Suppose that $f$ a strict concave function. In accordance with Lemma \ref{strictlyconcave}, 
\begin{align} \label{ASclass}
f(\bx+ w\bar \bp) < f(\bx) + w f'(\bx, \bar \bp) = f(\bx) + w\nabla f(\bx)^\top \bar\bp, 
\end{align}
which violates that $f$ satisfies (DLC) at $\bx, \bar \bp$ and $w \geq 0$. Hence, $f$ is not a strict concave function.

Assuming that $f$ is not a strict concave function: 
Suppose that $f$ does not satisfy (DLC). From Lemma \ref{strictlyconcave1}, $f$ is a strictly concave, which is a contradiction to $f$ is not a strictly concave. Therefore, $f$ must satisfy (DLC).

\qed

\section{Lemma \ref{dppmmonoprop}} \label{lemma-1}

From the definition of (DLC) and (\ref{updatedppm}),  we obtain
\begin{align}  \label{dppmmonoprop0}
[\nabla f(\bx) - \nabla f(\bu)]^\top (\bx - \bu)  & =  - w [\nabla f(\bx) - \nabla f(\bu)]^\top \bar \bp  \nonumber \\
& = [(\nabla f(\bx)^\top \bar \bp) - (\nabla f(\bu)^\top \bar \bp)]\bar \bp^\top (\bx - \bu) \geq 0.
\end{align} 
Without a loss of the generality, suppose that $w_2  > w_1$. (\ref{dppmmonoprop1}) is the result of  substituting $\bu_1$ and $\bu_2 =  \bu_1 + (w_2 - w_1) \bar \bp$ for $\bx$ and $\bu$ in  (\ref{dppmmonoprop0}), respectively.

$\bar \bp^\top \nabla f(\bx + w \bar \bp)$ is a monotonic function can be obtained by substituting $\bu_1 - \bu_2 = (w_1 - w_2) \bar \bp$ into (\ref{dppmmonoprop1}) to yield 
\begin{align*}
(w_1 - w_2)  (\nabla f(\bu_1) - \nabla f(\bu_2))^\top \bar \bp \geq 0.
\end{align*}
Since $w_2  > w_1$, 
\begin{align*}
\nabla f(\bu_1)^\top \bar \bp \leq  \nabla f(\bu_2)^\top \bar \bp.
\end{align*}

(\ref{dppmmonoprop1}) can be obtained from supposing that $\bar \bp^\top \nabla f(\bx + w \bar \bp)$ is an increasing function of $w \in [0, v(\bx, \bar \bp)]$, with
\begin{align*} 
[\bar \bp^\top \nabla f(\bu_1) - \bar \bp^\top \nabla f(\bu_2)] \bar \bp^\top (\bu_1 - \bu_2) = 
[\bar \bp^\top \nabla f(\bu_1) - \bar \bp^\top \nabla f(\bu_2)] (w_1 - w_2) \geq 0.
\end{align*}

\section{Lemma \ref{dppmlem}} \label{DPPMdescent}

We let $\bu= \bx- t (\nabla f(\bu)^\top \bar \bp) \bar \bp$ and use the convexity property for any $\bz \in (\bx, \bx + v(\bx) \bar \bp]$ to obtain
\begin{align*}
f(\bz) & \geq f(\bu) + \nabla f(\bu)^\top (\bz - \bu) \\
& = f(\bu)  + \nabla f(\bu)^\top  (\bz - \bx + t(\nabla f(\bu)^\top \bar \bp) \bar \bp).
\end{align*}
Re-writing the above yields 
\begin{align} \label{DPPMmain}
f(\bu) & \leq f(\bz) + \nabla f(\bu)^\top (\bx - \bz - t (\nabla f(\bu)^\top \bar \bp) \bar \bp \nonumber \\
& = f(\bz)  + \nabla f(\bu)^\top (\bx  - \bz) - t | \nabla f(\bu)^\top \bar \bp|^2.
\end{align}

Allowing $\bz = \bx$ obtains 
\begin{align} \label{dppmsufficient}
f(\bu) \leq f(\bx) - t | \nabla f(\bu)^\top \bar \bp|^2.
\end{align}
The conclusion follows by letting $\bu = \bx^{k+1}$, $\bx = \bx^k$, $t = t_k$, and $\bar \bp = \bar \bp_k$ in (\ref{dppmsufficient}).

\section{Theorem \ref{dppmlimitpoints}} \label{dppmlimitpointsproof}

(i) The entire sequence converges to a finite value of $f$ is a consequence of Lemma 4.1 \cite{zangwill1969nonlinear}, due to the fact that $\{\bx^k\}$ is a bounded sequence 
and $\{ f(\bx^k)\}$ is monotonically non-increasing that converges to a finite value of $f$. \\

(ii) To arrive at a contradiction, we assume that there exists a convergent subsequent to a non-critical point of $f$.  We first show that, under the assumption, we obtain $c_1 > 0$.
Let $\tilde K$ be the index set of a subsequence of iterates such that
\begin{align} \label{dppmx}
\{\bx^k\}_{k \in \tilde K, k \rightarrow \infty} \rightarrow \bx^*
\end{align}
and $\nabla f(\bx^*) \neq 0$.
From Lemma \ref{sdppmlem},  we obtain
\begin{align*}
\{c_{1, k}\}_{k \in \tilde K, k \rightarrow \infty} \rightarrow c_1.
\end{align*}
From (DLC) and $\nabla f(\bx^*) \neq 0$, there exists $\bar p_*$ such that $-\bar p_*^\top \nabla f(\bx^*) 
\geq - \bar p_*^\top \nabla f (\tilde \bx)$ where $\tilde \bx  = \bx^* + w^* \bar \bp_*$ and $w^* > 0$. From (\ref{dppmwolfe2}),  $c_1$ cannot be zero; otherwise, $- \bar p_*^\top \nabla f (\tilde \bx^*) = 0$. This violates the fact that $w^*$ is the solution to (\ref{dppmwvalue}).

From Lemma \ref{dppmlem} and $f$ is bounded from below, $\{ f(\bx^k)\}$ is monotonically non-increasing and converges to a finite value of $f$. Hence, 
\begin{align} \label{asymptotic}
f(\bx^k) - f(\bx^{k+1}) \rightarrow 0.
\end{align}
From Lemma \ref{sdppmlem} and the conclusion that $c_1 > 0$ under the assumption $\nabla f(\bx^*) \neq 0$, we obtain  
\begin{align}\label{dppmarmijo0}
f(\bx^k) - f(\bx^{k+1}) \geq - c_1 w_k \nabla f(\bx^k)^\top \bar \bp_k \;\; \forall k
\end{align}
Hence, from (\ref{asymptotic}), 
\begin{align} \label{dppmconvergence}
-w_k \nabla f(\bx^k)^\top \bar \bp_k \rightarrow 0.
\end{align}
$\{\bar \bp_k \}_{k \in \tilde K}$ is bounded that there is a subsequence $K$ of $\tilde K$ such that 
\begin{align} \label{dppmconvergence4}
\{\bar \bp_k\}_{k \in  K} \rightarrow \bar \bp.
\end{align}

$\{\bar \bp_k\}$ is gradient-related and therefore
\begin{align}\label{dppmconvergence5}
\lim\sup_{k \in K, k \rightarrow \infty} \nabla f(\bx^k)^\top \bar \bp_k < 0.
\end{align}
From (\ref{dppmconvergence5}) and (\ref{dppmconvergence}), we obtain
\begin{align} \label{dppmconvergence3}
\{ w_k \}_{k \in K} \rightarrow 0.
\end{align}

As shown in Figure \ref{mvt}, let $\alpha_k' > 0$ be the smallest intersecting value of $\alpha_k$ with
\begin{align} \label{alpha1}
f(\bx^k + \alpha_k' \bar \bp_k) = f(\bx^k) + \alpha_k' c_1 \nabla f(\bx^k)^\top \bar \bp_k.
\end{align}
The fact that $f \in \Gamma$ implies that there exists $\hat w_k \geq w_k$ in an interval around $\alpha_k'$, such that 
\begin{align*} 
f(\bx^k) - f(\bx^k + \hat w_k \bar \bp_k) < - c_1 \hat w_k \nabla f(\bx^k)^\top \bar \bp_k
\end{align*}
and 
\begin{align} \label{dppmmeanvalue}
\frac{f(\bx^k) - f(\bx^k +\hat w_k \bar \bp_k )}{\hat w_k} < -c_1 \nabla f(\bx^k)^\top \bar \bp_k.
\end{align}
Using the mean value theory, there exists  $\hat w_k' \in [0, \hat w_k]$ and (\ref{dppmmeanvalue}) can be expressed as 
\begin{align}\label{dppmconvergence2}
- \nabla f(\bx^k + \hat w_k' \bar \bp_k)^\top \bar \bp_k  < -c_1 \nabla f(\bx^k)^\top \bar \bp_k.
\end{align}
We can choose integer $l(k) \geq 1$ with $\hat w_k  = \frac{w_k}{l(k) \eta_k} \geq \alpha_k'$ with $\eta_k \in (0, 1)$ to satisfy (\ref{dppmconvergence2}). In accordance with (AS1), $l(k) \leq \frac{w_k}{\eta_k u_k} \leq \frac{1}{\eta_k} \leq \frac{1}{\eta}$ is bounded. From (\ref{dppmconvergence3}), we can have $\hat w_k \rightarrow 0$ when $k \in K$ and $k \rightarrow \infty$.  Notice that $\alpha_k'$ is defined as (56). By the mean value theorem, we know that there exists $\beta_k' \in [0, \alpha_k']$ such that
\begin{align} \label{meanvalue}
c_1 \nabla f(\bx^k)^\top \bar \bp_k = \phi'(\beta_k').
\end{align}

Following (\ref{meanvalue}) and $\lim_{k\rightarrow\infty}\alpha_k'  = 0$, we obtain
\[c_1\nabla f(\bx^*)^\top \bar \bp = \nabla f(\bx^*)^\top \bar \bp\]
where the last equality denotes the slope of $\phi(0)$ in direction $\bar \bp$.
Since $\lim\sup_{k \in K, k \rightarrow \infty} \nabla f(\bx^k)^\top \bar \bp_k < 0$, we obtain $c_1 = 1$. This is a contradiction to $c_1<1$ (Lemma \ref{sdppmlem}).

\qed

\section{Corollary \ref{dppmfejer}} \label{dppmfejerproof}

Theorem \ref{dppmlimitpoints} (i) allows us to assign $f(\bx^*) = f^*$ for any $x^* \in C(\bx^0)$. In accordance with the DPPM, we have
\begin{align} \label{DPPMdistance}
\| \bu - \bx^* \|^2 & = \| \bx - t (\nabla f(\bu)^\top \bar \bp) \bar \bp- \bx^* \|^2  \nonumber \\
& = \|\bx - \bx^* \|^2 - 2t (\nabla f(\bu)^\top \bar \bp) \bar \bp^\top (\bx - \bx^*) + t^2 |\nabla f(\bu)^\top \bar \bp |^2.
\end{align}
We let $\bx$ denote an iterator $\bx^k$ with $k \geq k_0$ and let $\bu= \bx- t (\nabla f(\bu)^\top \bar \bp) \bar \bp$. From the fact that $\bx^k \in \mathcal D$ for $k \geq k_0$ and $f$ over $\mathcal D$ is convex, 
we obtain, for $\bz \in \mathcal D$ 
\begin{align}
f(\bz) \geq f(\bu) + \nabla f(\bu)^\top (\bz - \bu) 
 = f(\bu)  + \nabla f(\bu)^\top  (\bz - \bx + t(\nabla f(\bu)^\top \bar \bp) \bar \bp).
\end{align}
This equation can be re-written as  
\begin{align} \label{DPPMmain0}
f(\bu)  \leq  f(\bz)  + \nabla f(\bu)^\top (\bx  - \bz) - t | \nabla f(\bu)^\top \bar \bp|^2.
\end{align}
If (\ref{gradientdecomp}) is applied to $\nabla f(\bu)^\top (\bx  - \bz)$, then we have 
\begin{align} \label{dppmcrossing}
\nabla f(\bu)^\top (\bx - \bz) & = \| \nabla f(\bu)\| \overline{\nabla f(\bu)}^\top (\bx - \bz) 
 = \| \nabla f(\bu)\| (\alpha \bar \bp + \beta \bv)^\top (\bx - \bz)  
\end{align}
Substituting (\ref{dppmcrossing}) into (\ref{DPPMmain0}) yields
\begin{align} \label{fejer1}
f(\bu) & \leq f(\bz)  + \nabla f(\bu)^\top (\bx  - \bz) - t | \nabla f(\bu)^\top \bar \bp|^2 \nonumber \\
& = f(\bz) + \| \nabla f(\bu)\| (\alpha \bar \bp + \beta \bv)^\top (\bx - \bz)   - t | \nabla f(\bu)^\top \bar \bp|^2.
\end{align}
From (\ref{gradientdecomp}), we can obtain
\begin{align}
\alpha = \bar \bp^\top \overline{\nabla f(\bu)} - \beta \bar \bp^\top \bv.
\end{align}
Substituting this $\alpha$ into (\ref{fejer1}) yields
\begin{align} \label{fejer11}
f(\bu)  & \leq f(\bz) + \| \nabla f(\bu)\| ((\bar \bp^\top \overline{\nabla f(\bu)} - \beta \bar \bp^\top \bv) \bar \bp + \beta \bv)^\top (\bx - \bz)   - t | \nabla f(\bu)^\top \bar \bp|^2 \nonumber \\
& =  f(\bz) + (\bar \bp^\top \nabla f(\bu) - \beta  \| \nabla f(\bu)\|\bar \bp^\top \bv) \bar \bp + \beta  \| \nabla f(\bu)\|\bv)^\top (\bx - \bz)   - t | \nabla f(\bu)^\top \bar \bp|^2 
\end{align}
We can now substitute $\bx^* \in \mathcal D$ for $\bz$ in (\ref{fejer11}) to obtain  
\begin{align} \label{fejer2}
f(\bu) & \leq  f^* +  (\nabla f(\bu)^\top \bar \bp) \bar \bp^\top + \beta \| \nabla f(\bu)\| \bv^\top - (\beta  \| \nabla f(\bu)\| \bar \bp^\top \bar \bv) \bar \bp^\top)(\bx - \bx^*) - t | \nabla f(\bu)^\top \bar \bp|^2.
\end{align}
Multiplying both sides of (\ref{fejer2}) by $2t$, and replacing $ 2t (\nabla f(\bu)^\top \bar \bp) \bar \bp^\top (\bx - \bx^*)$ in accordance with (\ref{DPPMdistance}), we obtain 
\begin{align} \label{DPPMvalue0}
2t (f(\bu) - f^*) & \leq  \|\bx - \bx^* \|^2 - \| \bu - \bx^* \|^2 - t^2 |\nabla f(\bu)^\top \bar \bp |^2 \nonumber \\
& + 2 t\beta \| \nabla f(\bu)\|\bv^\top (\bx - \bx^*) - 2t \beta \| \nabla f(\bu)\| (\bar \bp^\top \bv) \bar \bp^\top (\bx - \bx^*)
\end{align}
Note that $\bx - \bx^* = -\bar \bp^* \|\bx - \bx^*\|$, 
\begin{align} \label{DPPMvalue4}
2t \beta  \| \nabla f(\bu)\| (\bv^\top  - (\bar \bp^\top \bv) \bar \bp^\top) (\bx - \bx^*) &  = - 2 \|\bx - \bx^*\|t \beta  \| \nabla f(\bu)\| (\bv^\top  - (\bar \bp^\top \bv) \bar \bp^\top)\bar \bp^* \nonumber \\
 & = - 2 \|\bx - \bx^*\|t \beta  \| \nabla f(\bu)\| (\bv^\top \bar \bp^* - (\bar \bp^\top \bv)^2 ) \nonumber \\
 & =  2 \|\bx - \bx^*\|t \beta  \| \nabla f(\bu)\| (-\| \bv \|^2 + (\bar \bp^\top \bv)^2 ) 
\end{align}
Because $-\| \bv \|^2 + (\bar \bp^\top \bv)^2  \leq 0$, if $\beta \geq 0$, then (\ref{DPPMvalue4}) is smaller than zero; hence, (\ref{DPPMvalue0}) becomes
\begin{align} \label{DPPMvalue2}
2t (f(\bu) - f^*) & \leq  \|\bx - \bx^* \|^2 - \| \bu - \bx^* \|^2 - t^2 |\nabla f(\bu)^\top \bar \bp |^2.
\end{align}
From the fact that the DPPM chooses $\cos (\bar \bp, \bar \bv) > 0$,  we obtain $\beta \geq 0$, in accordance with (\ref{cosine}). 
Since $w = t |\nabla f(\bu)^\top \bar \bp|$, (\ref{DPPMvalue2}) can be expressed as 
\begin{align} \label{DPPMvalue1}
0 \leq 2t (f(\bu) - f^* ) + w^2 \leq \|\bx - \bx^* \|^2 - \| \bu - \bx^* \|^2. 
\end{align}
Hence, 
\begin{align*}
\| \bu- \bx^*\| \leq \|\bx - \bx^*\|.
\end{align*}
Therefore,  we obtain
\begin{align}  \label{fejer3}
\| \bx^{k+1}- \bx^*\| \leq \|\bx^k - \bx^*\|,
\end{align}
for any $k \geq k_0$. Since $\bx^*$ is any critical point in $C(\bx^0)$. We obtain the conclusion.

\qed

%

\section{DPPM Acceleration} \label{sec:acc}

For the sake of convenience, equations  (\ref{dppmdescent})  and (\ref{DPPMvalue1}) are re-stated, respectively:
\begin{align}  \label{dppmaccer1}
f(\bu) \leq f(\bx) - t | \nabla f(\bu)^\top \bar \bp|^2, 
\end{align}
\begin{align}  \label{dppmaccer0}
t( f(\bu) -f^*)   \leq \|\bx - \bx^* \|^2 - \| \bu - \bx^* \|^2.
\end{align} 

Let $z =  | \nabla f(\bu)^\top \bar \bp|$. Multiply (\ref{dppmaccer1}) by $t(1-\theta)$ where $\theta \in [0, 1]$ and multiply (\ref{dppmaccer0}) by $\theta$ and then add the resultants together yields
\begin{align*}
t f(\bu) & \leq (1-\theta) t f(\bx) + t \theta f^* - t^2 (1-\theta) z^2 + \theta (\|\bx - \bx^* \|^2 - \| \bu - \bx^* \|^2) \\
& \leq (1-\theta) t f(\bx) + t \theta f^* - t^2 (1-\theta) z^2 - 2 tz \theta \bar \bp^\top (\bx - \bx^*) - \theta t^2 z^2.
\end{align*}
The last equality is derived with $\bu= \bx - t (\nabla f(\bu)^\top \bar \bp) \bar \bp = \bx + t z \bar \bp$ in accordance with $\nabla f(\bu)^\top \bar \bp \leq 0$.
Following the inequality, we obtain
\begin{align*}
f(\bu) & \leq (1- \theta) f(\bx) + \theta f^* - 2\theta z \bar \bp^\top (\bx - \bx^*) - t z^2 \\
& = (1- \theta) f(\bx) + \theta f^* + \frac{\theta^2}{t} (\|\bx - \bx^* \|^2 - \| \bx + \frac{tz \bar \bp}{\theta} - \bx^* \|^2).
\end{align*}
Hence, 
\begin{align} \label{dppmaccer2}
f(\bu) - f^* \leq (1-\theta) f(\bx) - (1-\theta) f^* + \frac{\theta^2}{t}  (\|\bx - \bx^* \|^2 - \| \bx + \frac{tz \bar \bp}{\theta} - \bx^* \|^2).
\end{align}
Letting $\bv = \bx + \frac{tz \bar \bp}{\theta}$ and recall that $
\bu = \bx + tz \bar \bp$,  we obtain  
\begin{align} \label{vextra}
\bv = \bx + \frac{1}{\theta} (\bu - \bx) = (1 - \frac{1}{\theta}) \bx + \frac{1}{\theta} \bu.
\end{align}

Dividing both sides of (\ref{dppmaccer2}) by $\theta^2$, and re-arranging the terms yields
\begin{align} \label{dppmaccer2}
\frac{1}{\theta^2} (f(\bu) - f^*) + \frac{1}{t} \|\bv - \bx^*\|^2 \leq \frac{1-\theta}{\theta^2}[f(\bx) - f^*] + \frac{1}{t} \|\bx - \bx^*\|^2.
\end{align}
If we retrieve the iteration index of the above and replace $\bu, \bv, \bx, \theta$ with $\bx^k$, $\bv^k$, $\bv^{k-1}$, and $\theta_k$, respectively, Equation (\ref{dppmaccer2}) can be expressed as
\begin{align} \label{dppmaccer3}
\frac{1}{\theta_k^2} (f(\bx^k) - f^*) + \frac{1}{t} \|\bv^k - \bx^*\|^2 \leq \frac{1-\theta_k}{\theta_k^2}[f(\bv^{k-1}) - f^*] + \frac{1}{t} \|\bv^{k-1} - \bx^*\|^2.
\end{align}
Letting 
\begin{align}
\theta_i = \frac{2}{i+1},
\end{align}
we have the desired properties: $\theta_1 = 1$, $\theta_i \in (0, 1)$, and 
\begin{align} \label{dppmthetacon}
\frac{1-\theta_i}{\theta_i^2} \leq \frac{1}{\theta_{i-1}^2}.
\end{align}
Repeatedly using (\ref{dppmaccer3}) and (\ref{dppmthetacon}), we can obtain
\begin{align} \label{dppmaccer4}
\frac{1}{\theta_k^2} (f(\bx^k) - f^*) + \frac{1}{t} \|\bv^k - \bx^*\|^2 & \leq \frac{1}{\theta_{k-1}^2}[f(\bv^{k-1}) - f^*] + \frac{1}{t} \|\bv^{k-1} - \bx^*\|^2 \nonumber \\
& \leq \frac{1-\theta_{k-1}}{\theta_{k-1}^2}[f(\bv^{k-2}) - f^*] + \frac{1}{t} \|\bv^{k-2} - \bx^*\|^2 \nonumber\\
& \leq \frac{1}{\theta_{k-2}^2}[f(\bv^{k-2}) - f^*] + \frac{1}{t} \|\bv^{k-2} - \bx^*\|^2 \nonumber\\
& \vdots \nonumber\\
& \leq \frac{1}{\theta_{1}^2}[f(\bv^{1}) - f^*] + \frac{1}{t} \|\bv^{1} - \bx^*\|^2 \nonumber\\
& \leq \frac{1- \theta_1}{\theta_{1}^2}[f(\bv^{0}) - f^*] + \frac{1}{t} \|\bv^{0} - \bx^*\|^2 \nonumber \\
& \leq \frac{1}{t} \|\bv^{0} - \bx^*\|^2.
\end{align}
Deduced from (\ref{dppmaccer4}) and $\bv^0 = \bx^0$,  
\begin{align} \label{dppmaccer5}
f(\bx^k) - f^* \leq \frac{\theta_k^2}{t} \|\bx^{0} - \bx^*\|^2.
\end{align}
The order of $k$ to achieve $\epsilon$-suboptimal solution of $f^*$ is $\frac{1}{\sqrt{\epsilon}}$.

\qed

%
%
%

Eq. (\ref{vextra}) indicates that $\bv^k$ is an extrapolation of $\bx^k$ and $\bu^k$ with parameter $\theta_k = \frac{2}{k+1}$. In the acceleration, $\bu^k$ is the update of $\bx^k$ using the DPPM with constant $t$ and the next iteration $\bx^{k+1}$ the the DPPM is the extrapolation output $\bv^k$. The process repeats with initial setting $\bx^0 = \bv^0$ and constant $t$.

%

\section{Proof of Lemma \ref{PROPOSITION_1}} \label{inversequad}

(i) Let $\{\bp_1= \bp, \bp_2 , \cdots \bp_n \}$ 
be an orthonormal basis of $\mathbb{R}^{n}$. If $(1+t \bar{\bp}^{\top} Q \bar{\bp} )\ne 0$, then 
\begin{equation}
[\bI + t\bar{\bp} \bar{\bp}^{\top} Q]\bp_i
= \bp_i + (t{\bp_1}^{\top} Q \bp_i) \bp_1,
\end{equation}
for each $i \in \{1,2,\cdots , n\}$. Hence, $\bI + t\bar{\bp} \bar{\bp}^{\top} Q$ is invertible. \\
Suppose $[\bI + t\bar{\bp} \bar{\bp}^{\top} Q]^{-1} =\bI+M$ 
for some $M \in \mathbb{R}^{n\times n}$ and suppose 
\begin{eqnarray}\label{Def_AB}
A =  t\bar{\bp}^{\top} Q \text{ and }
B =  \bar{\bp}^{\top} M. \nonumber
\end{eqnarray}
Then, we obtain
\begin{equation}\label{Prop_pptM}
\bI+M = \bI+ \bar{\bp} B +(\bI - \bar{\bp} \bar{\bp}^{\top})M.
\end{equation}
From the fact that 
$[\bI + t\bar{\bp} \bar{\bp}^{\top} Q][\bI+M] = [\bI+\bar{\bp}A][\bI+M] = \bI $ and (\ref{Prop_pptM}), we obtain
\begin{eqnarray}\label{Prop_B}
0 &=& [\bI+\bar{\bp}A][\bI+M] -\bI \nonumber\\
&=& M +\bar{\bp}A  [\bI+  M]\nonumber\\
&=& (\bI - \bar{\bp} \bar{\bp}^{\top}) [ M +\bar{\bp}A  [\bI+  M]]\nonumber\\
&=& (\bI - \bar{\bp} \bar{\bp}^{\top}) M.
\end{eqnarray}
In accordance with (\ref{Prop_pptM}) and (\ref{Prop_B}), we have
\begin{eqnarray}\label{Prop_pB}
0 &=& [\bI+\bar{\bp}A][\bI+M] -\bI \nonumber\\
&=& [\bI+\bar{\bp}A][\bI+\bar{\bp} B ] -\bI \nonumber\\
&=& \bar{\bp}B+\bar{\bp}A +  \bar{\bp}A \bar{\bp}B \nonumber\\
&=& \bar{\bp}B+ t\bar{\bp} \bar{\bp}^{\top} Q +   t\bar{\bp} \bar{\bp}^{\top} Q \bar{\bp}B \nonumber\\
&=&  t\bar{\bp} \bar{\bp}^{\top} Q + (1+t \bar{\bp}^{\top} Q \bar{\bp})\bar{\bp}B. 
\end{eqnarray}
From (\ref{Prop_pB}), we obtain $\bar{\bp}B = - \frac{t}{(1+t \bar{\bp}^{\top} Q \bar{\bp} )}\bar{\bp} \bar{\bp}^{\top} Q $ and, hence, 
\begin{equation}
[\bI + t\bar{\bp} \bar{\bp}^{\top} Q]^{-1} 
= \bI+M = \bI+\bar{\bp}B  
=\bI - \frac{t}{(1+t \bar{\bp}^{\top} Q \bar{\bp} )}\bar{\bp} \bar{\bp}^{\top} Q.
\end{equation}

(ii) Let $\bar{\bp} \perp \bar \bp_{\perp}$ and let $a \bar{\bp} + b \bar \bp_{\perp}$ ($a^2 + b^2 = 1$) be an eigenvector of  
$\bI + t\bar{\bp} \bar{\bp}^{\top} Q$ with corresponding eigenvalue $\lambda$. Then, 
\begin{equation}\label{Prop_pperp}
[\bI + t\bar{\bp} \bar{\bp}^{\top} Q](a \bar{\bp} + b \bar \bp_{\perp})
= (a+at\bar{\bp}^{\top} Q \bar{\bp} + bt\bar{\bp}^{\top} Q\bar \bp_{\perp} )\bar{\bp} 
+ b \bar \bp_{\perp}
= \lambda (a \bar{\bp} + b \bar \bp_{\perp}).
\end{equation}
\noindent
From this equation, we have
\begin{itemize}
\item
If $b \ne 0$, then $\lambda = 1$.
\item
If $b = 0$ (a.k.a $a = 1$), then 
$\lambda = 1+t \bar{\bp}^{\top} Q \bar{\bp}$ and eigenvector is $\bar \bp$.
\end{itemize}
\noindent
Hence, $\bar \bp$ is the eigenvector of $[\bI + t\bar{\bp} \bar{\bp}^{\top} Q]$ with corresponding eigenvalue ${1+ t \bar \bp^\top Q \bar \bp}$ and the other eigenvectors of $[\bI + t\bar{\bp} \bar{\bp}^{\top} Q]$ have corresponding eigenvalues equal to $1$.

\end{document}